\documentclass[a4paper,10pt]{article}
\usepackage{amssymb}
\newtheorem{Th}{Theorem}
\newtheorem{Prop}{Proposition}

\newtheorem{Lm}{Lemma}

\newtheorem{Dfi}{Definition}
\newtheorem{Rm}{Remark}

\newcommand{\be}{\begin{equation}}
\newcommand{\ee}{\end{equation}}
\newcommand{\R}{\mathbb{R}}
\newcommand{\N}{\mathbb{N}}

\newcommand\res{\mathop{\hbox{\vrule height 7pt width .5pt depth 0pt
\vrule height .5pt width 6pt depth 0pt}}\nolimits}

\newcommand{\reset}{\setcounter{equation}{0}\setcounter{Th}{0}\setcounter{Prop}{0}\setcounter{Co}{0}
\setcounter{Lm}{0}\setcounter{Rm}{0}}

\def\ti{\tilde}
\def\lf{\left}
\def\rg{\right}

\def\al{\alpha}
\def\la{\lambda}

\def\ep{\varepsilon}
\def\ds{\displaystyle}
\def\ov{\overline}

\def\om{\omega}
\def\p{\partial}

\def\res{\mathop{\hbox{\vrule height 7pt width .5pt 
depth 0pt\vrule height .5pt width 6pt depth 0pt}}\nolimits}
\begin{document}
\title{Sequences of Smooth Global Isothermic Immersions.}
\author{Tristan Rivi\`ere\footnote{Department of Mathematics, ETH Zentrum,
CH-8093 Z\"urich, Switzerland.}}
\date{ }
\maketitle

{\bf Abstract :} {\it In the present work we study the behavior of sequences of smooth global isothermic immersions of a given
 closed surface and having a
uniformly bounded total curvature. We prove that, if the conformal class of this sequence is bounded in the
Moduli space of the surface, it weakly converges in $W^{2,2}$ away from finitely many points, modulo extraction of a subsequence, to a possibly branched weak isothermic immersion of this surface.
Moreover, if this limit happens to be smooth away from the branched points, we give an optimal description of the possible loss of strong compactness  of such a subsequence by proving that, beside possibly finitely many
atomic concentrations, the defect measure associated to the $L^2$ norm
of the second fundamental form is ''transported'' along exceptional directions given by some holomorphic quadratic forms associated the limiting surface. We give examples where such a loss of compactness, invariant along such exceptional directions, eventually happen.}

\medskip

\noindent{\bf Math. Class.} 35L51, 35L65, 35R01, 30C70, 53A30, 58E30, 49Q10, 35J35, 35J48, 35J50.

\section{Introduction to Global Isothermic Immersions.}

\subsection{The origin of isothermic in the XIXth century's surface geometry in ${\R}^3$ and its generalization to arbitrary codimensions.}

The notion of isothermic surfaces has been introduced in the second half of the XIX century and was in particular studied
by E. Bour, E.B. Christoffel and G. Darboux in the context of conjugated famillies of surfaces.  The issue was to find pairs of distinct, non homothetic, immersions into ${\R}^3$,
 $\vec{\Phi}$ and $\vec{L}$ of the 2 dimensional disc $D^2$  ''dual'' to each other in the following sense\footnote{Darboux formulated the problem this way (see \cite{Da2}) :  {\it Proposons nous de rechercher tous les cas dans lesquels la correspondance par plan tangents parall\`eles \'etablie entre deux surfaces peut donner une repr\'esentation conforme ou un trac\'e g\'eographique de l'une des surfaces sur l'autre.} } :
 \be
 \label{I.1}
 \p_{x_i}\vec{\Phi}\quad\mbox{ is parallel to }\quad\p_{x_i}\vec{L}\quad\quad\mbox{ for } i=1,2
 \ee
 and the two induced metric on $D^2$ are conformal to each other :
 \be
 \label{I.2}
 \vec{L}^\ast g_{{\R}^3}=e^{2u}\ \vec{\Phi}^\ast g_{{\R}^3}
 \ee
 where $g_{{\R}^3}$ is the standard metric on ${\R}^3$ and $u$ is an arbitrary function on $D^2$.
 
 \medskip
 
E. Bour and   E.B. Christoffel proved respectively in \cite{Bou} and \cite{Chr} that the non trivial solutions to this question are immersions which posses around every point conformal (or isothermic) coordinates such that the coordinate directions are principal (or curvature lines). In other words if $\vec{n}_{\vec{\Phi}}$ denotes the 
 Gauss Map of such an immersion
 \[
 \vec{n}_{\vec{\Phi}}:=\frac{\p_{x_1}\vec{\Phi}\times\p_{x_2}\vec{\Phi}}{|\p_{x_1}\vec{\Phi}\times\p_{x_2}\vec{\Phi}|}
 \]
 around each point there exists $(x_1,x_2)$ coordinates such that the induced metric is conformal
 \be
 \label{I.3}
 \vec{\Phi}^\ast g_{{\R}^3}=e^{2\la}\ \lf[ dx_1^2+dx_2^2\rg]
 \ee
 and 
 \be
 \label{I.4}
 <\p_{x_1}\vec{n}_{\vec{\Phi}},\p_{x_2}\vec{\Phi}>= <\p_{x_2}\vec{n}_{\vec{\Phi}},\p_{x_1}\vec{\Phi}>=0
\ee
where $<\cdot,\cdot>$ denotes the scalar product in ${\R}^3$, which also means that the second fundamental form is diagonal in these conformal coordinates :
\[
\vec{{\mathbb I}}= -e^{-2\la}\ \lf[<\p_{x_1}\vec{n}_{\vec{\Phi}},\p_{x_1}\vec{\Phi}>\ dx_1^2+ <\p_{x_2}\vec{n}_{\vec{\Phi}},\p_{x_2}\vec{\Phi}>\ dx_2^2\rg]\ \vec{n}_{\vec{\Phi}}\quad.
\]
where $e^\la=|\p_{x_1}\vec{\Phi}|=|\p_{x_2}\vec{\Phi}|$. If (\ref{I.3}) and (\ref{I.4}) hold one says that the curvature lines are {\it isothermic} and, following Darboux,  such a surface is called {\it isothermic surface}. Since that time example of
isothermic surfaces were known such as {\it axially symmetric surfaces} or {\it constant mean curvature surfaces} including of course {\it minimal surfaces}.

\medskip

In order to extend the notion of isothermic surfaces to immersions into ${\R}^n$ for an arbitrary $n>2$ we need to reformulate  the pair of constraints (\ref{I.1}) and (\ref{I.2})
or equivalently the pair of constraints (\ref{I.3}) and (\ref{I.4}) but also to relax slightly this assumption.

\medskip

We recall the definition of the {\it Weingarten form} $\vec{h}_0$ of an immersion $\vec{\Phi}$ into ${\R}^3$, in an arbitrary choice of complex coordinates,
\[
\begin{array}{l}
\vec{h}_0:= -e^{-2\la}\ <\p_z\vec{n}_{\Phi},\p_z\vec{\Phi}>\ dz\otimes dz\\[5mm]
\ds\quad\quad=-\frac{e^{-2\la}}{4}\lf[<\p_{x_1}\vec{n}_{\vec{\Phi}},\p_{x_1}\vec{\Phi}>-<\p_{x_2}\vec{n}_{\vec{\Phi}},\p_{x_2}\vec{\Phi}>-2\, i\ <\p_{x_1}\vec{n}_{\vec{\Phi}},\p_{x_2}\vec{\Phi}>\rg]\ \vec{n}_{\vec{\Phi}}\ {dz\otimes dz}
\end{array}
\]
where $z=x_1+i x_2$ and $\p_z:={2^{-1}}\ [\p_{x_1}-i\p_{x_2}]$.

\medskip

\noindent As observed in \cite{Ri3} we have the following result.

\medskip

\begin{Prop}
\label{pr-I.1} A conformal immersion $\vec{\Phi}$ of the disc $D^2$ into ${\R}^3$ satisfies, around each point, except possibly a discrete subset of $D^2$,  (\ref{I.4})  in some other local conformal
chart $(y_1,y_2)$  if and only if there exists a non zero holomorphic function $f(z)$ on $D^2$ such that
\be
\label{I.6}
\Im\lf(\ov{f(z)}\ \vec{H}_0\rg)=0\quad,
\ee
where $\vec{H}_0:=-4^{-1}\,e^{-2\la}\ \lf[<\p_{x_1}\vec{n}_{\vec{\Phi}},\p_{x_1}\vec{\Phi}>-<\p_{x_2}\vec{n}_{\vec{\Phi}},\p_{x_2}\vec{\Phi}>-2\, i\ <\p_{x_1}\vec{n}_{\vec{\Phi}},\p_{x_2}\vec{\Phi}>\rg]\ \vec{n}_{\vec{\Phi}}$ is the expression of $\vec{h}_0$ in the given conformal parametrization $\vec{\Phi}$ on the disc $D^2$.
\hfill $\Box$
\end{Prop}
Indeed, while changing conformal coordinates and taking $w(z):=y_1(z)+i y_2(z)$ the expression of $\vec{h}_0$ in these new coordinates becomes
\be
\label{I.7}
\vec{H}_0'\circ w= |\p_zw|^2\,(\p_z w)^{-2}\ \vec{H}_0\quad.
\ee
Away from the zeros of $f$, taking $w(z):=\sqrt{f(z)}$, (\ref{I.6}) becomes
\[
\Im(\vec{H}_0')=0\quad,
\]
which is exactly (\ref{I.4}).

\medskip

We introduce on the space $\wedge^{1-0}D^2\otimes\wedge^{1-0} D^2$ of $1-0\otimes1-0$ form on $D^2$ the following hermitian product\footnote{ This hermitian product integrated on $D^2$ is the {\it Weil Peterson product}.} depending on the conformal immersion $\vec{\Phi}$
\[
(\psi_1\ dz\otimes dz,\psi_2\ dz\otimes dz)_{WP}:= e^{-4\la}\ \ov{\psi_1(z)}\ \psi_2(z)
\]
where $e^\la:=|\p_{x_1}\vec{\Phi}|=|\p_{x_2}\vec{\Phi}|$. We observe that for a conformal change of coordinate $w(z)$ (i.e. $w$ is holomorphic in $z$) and for $\psi_i'$ satisfying
\[
\psi_i'\circ w\ dw\otimes dw=\psi_i\ dz\otimes dz
\]
one has, using the conformal immersion $\vec{\Phi}\circ w$ in the l.h.s.
\[
(\psi_1'\ dw\otimes dw,\psi_2'\ dw\otimes dw)_{WP}=(\psi_1\ dz\otimes dz,\psi_2\ dz\otimes dz)_{WP} 
\]
Using this change of coordinate rule, (\ref{I.6}) is equivalent to the following intrinsic characterization :
there exists an holomorphic section\footnote{ In complex coordinates $q=f(z)\ dz\otimes dz$ where $f$ is holomorphic and $q$ is called an holomorphic quadratic form.} $q$ of the bundle $\wedge^{1-0}D^2\otimes\wedge^{1-0} D^2$   such that
\be
\label{I.8}
\Im(q,\vec{h}_0)_{WP}=0\quad.
\ee

In codimension larger than 1 principal directions are not defined anymore and the XIXth century definition of isothermic immersions into ${\R}^3$
cannot be extended in a straightforward way for immersions into ${\R}^m$ ($m>3$).  However, considering a smooth immersions $\vec{\Phi}$ of an arbitrary 2-dimensional manifold $\Sigma$
into ${\R}^m$ one can still produce the global {\it Weingarten form} using local conformal charts as being the following global section of ${\R}^m\otimes\wedge^{1-0}\Sigma\otimes\wedge^{1-0} \Sigma$  :
\be
\label{I.9}
\begin{array}{l}
\vec{h}_0:= 2\,e^{-2\la}\,\pi_{\vec{n}}(\p^2_{z^2}\vec{\Phi})\ dz\otimes dz\\[5mm]
\ds\quad=\frac{e^{-2\la}}{2}\,\pi_{\vec{n}}\lf(\p^2_{x^2_1}\vec{\Phi}-\p^2_{x^2_2}\vec{\Phi}-2\, i\ \p^2_{x_1x_2}\vec{\Phi}\rg)\ \ {dz\otimes dz}
\end{array}
\ee
where $\pi_{\vec{n}}$ is the orthogonal projection onto the plane orthogonal to $\vec{\Phi}_\ast T\Sigma$.
%where the Gauss map $\vec{n}_{\Phi}$ is the unit simple 2-vector of $\wedge^2{\R}^m$ giving the oriented tangent space $\vec{\Phi}_\ast T\Sigma$ and $\res$ is the standard contraction operator between  2-vectors and vectors given by
%\[
%\forall\,\vec{a}\in\wedge^2{\R}^m\ ,\quad\forall\,\vec{b}\, ,\vec{c}\in\wedge^1{\R}^m\quad\quad<\vec{a}\res\vec{b},\vec{c}>=<\vec{a},\vec{b}\wedge\vec{c}>\quad.
%\] 
We can now introduce the natural generalization of global smooth isothermic surfaces into arbitrary euclidian space ${\R}^m$.
\begin{Dfi}
\label{df-I.1}
Let $\vec{\Phi}$ be a smooth immersion of a two dimensional manifold $\Sigma^2$ into ${\R}^m$. One says that $\vec{\Phi}$ is {\bf global isothermic} if
there exists an holomorphic quadratic form $q$ of the riemann surface issued from $\Sigma^2$ equipped with the pull back metric $g:=\vec{\Phi}^\ast g_{{\R}^m}$ 
of the standard metric $g_{{\R}^m}$ of ${\R}^m$ such that
\be
\label{I.10}
\Im(q,\vec{h}_0)_{WP}=0\quad.
\ee
where $\vec{h}_0$ is the Weingarten form of the immersion $\vec{\Phi}$ given by (\ref{I.9}).\hfill $\Box$
\end{Dfi}

\subsection{The role of Isothermic surfaces in the calculs of variations of the Willmore Lagrangian.}

In this work we are interested with analysis properties of Smooth global isothermic immersions. One of the main reasons why looking at the analysis of global isothemic
immersions comes from the fact that they may arise as degenerate critical point to the conformal constrained Willmore problem, as it has been shown in \cite{Ri3}.
In his 3 volumes book on differential geometry published by Springer around 1929 Wilhelm Blaschke, (see in particular the third volume
\cite{Bla}) proposed a theory
merging {\it minimal surface theory} and {\it conformal invariance}. This theory consists in studying the variations of the now so called
{\it Willmore Lagrangian} for surfaces. This lagrangian, $W$, is given by the $L^2$ norm of the mean curvature vector $\vec{H}$ of  an arbitrary immersion $\vec{\Phi}$ 
into the euclidian space ${\R}^m$ ($m\ge 3$) of a given
2-dimensional abstract manifold $\Sigma$  and integrated with respect to the induced metric\footnote{ $g:=\vec{\Phi}^\ast g_{{\R}^m}$ where $g_{{\R}^m}$
is the canonical flat metric of ${\R}^m$.} $g$
\be
\label{I.10a}
W(\vec{\Phi}):=\int_{\Sigma}|\vec{H}|^2\ dvol_g\quad.
\ee
Immersions satisfying $W(\vec{\Phi})<+\infty$ are called {\it immersions of finite total curvature}.

Minimal immersions, satisfying $\vec{H}\equiv 0$, are clearly critical points to $W$. Blaschke observed\footnote{This invariance was proved by Wilhelm Blaschke for $m=3$
and later on generalized by Bang-Yen Chen to arbitrary $m$} moreover the following conformal invariance of the
lagrangian $W$ : for any conformal diifeomorphism $\Psi$ of ${\R}^m\cup\{\infty\}$ into itself which does not send any point of $\Phi(\Sigma)$ to infinity one has
\be
\label{I.11}
W(\Psi\circ\vec{\Phi})=W(\vec{\Phi})\quad.
\ee
Hence, as a consequence, any composition of a minimal surface with a conformal diffeomorphism is still a critical point of $W$ without being necessarily minimal
anymore. Though the space of critical points of $W$ happens to be much broader than such compositions, Blaschke decided nevertheless to call such an immersion a {\it conformal minimal immersion}\footnote{ Probably in order to insist on the merging of the two requirements for this theory to include minimal surfaces and conformal invariance }. {\it Conformal minimal immersions} are nowadays
known under the denomination {\it Willmore surfaces}. Example of such surfaces are given for instance by minimal surfaces in ${\R}^m$ or stereographic projections
into ${\R}^m$ of minimal surfaces in $S^m$, constant mean curvature surfaces in ${\R}^3$  and the compositions of all these surfaces with conformal transformations .
It has been proven in \cite{Ri2} that an immersion $\vec{\Phi}$ is a critical point to  $W$ if and only if it satisfies
\be
\label{I.12}
d^{\ast_g}\lf[ d\vec{H}-3\, D\vec{H}+\star(\ast_{g}d\vec{n}_{\vec{\Phi}}\wedge\vec{H})\rg]=0
\ee
where $\ast_g$ is the Hodge operator on $\Sigma$ associated to the induced metric $g:=\vec{\Phi}^\ast g_{{\R}^m}$, $D\vec{H}$ is the covariant differentiation
of the section $\vec{H}$ of the normal bundle $(\vec{\Phi}_\ast T\Sigma)^\perp)$, it is also given by
\[
D\vec{H}:=\pi_{\vec{n}}(d\vec{H})
\]
where $\pi_{\vec{n}}$ denotes the orthogonal projection onto the fibers of $(\vec{\Phi}_\ast T\Sigma)^\perp$. Finally $\star$ denotes the Hodge operator from $\wedge^p{\R}^m$
into $\wedge^{m-p}{\R}^m$ for the canonical metric of ${\R}^m$ satisfying
\[
\forall \al\,\beta\in \wedge^p{\R}^m\quad\quad\al\wedge \star\beta=(\al,\beta)\ \ep_1\wedge\cdots\wedge\ep_m
\]
where $\ep_i$ is the canonical basis of ${\R}^m$ and $(\cdot,\cdot)$ denotes the canonical scalar product on $\wedge^p{\R}^m$.
In conformal coordinates for the induced metric $g$ equation (\ref{I.12}) becomes.
\be
\label{I.13}
div\lf(\nabla\vec{H}-3\,\pi_{\vec{n}}(\nabla\vec{H})+\star(\nabla^\perp\vec{n}_{\vec{\Phi}}\wedge\vec{H})\rg)=0\quad.
\ee

While exploring the existence and properties of critical points to the Willmore energy (\ref{I.11}) , or in other words while proceeding to the calculus of variation
of the Lagrangian $W$, it is natural to raise the question of the conformal class 
such an immersion defines on the abstract 2-manifold $\Sigma$. As a channel of  consequences it is then natural to explore minimizers or critical points to $W$
when the conformal class induced by $\vec{\Phi}^\ast g_{{\R}^m}$ is fixed. Assuming such a critical point is a non degenerate point for the conformal class constraint,
it has been proved in \cite{Ri3} that $\vec{\Phi}$ satisfies this time
\be
\label{I.14}
d^{\ast_g}\lf[ d\vec{H}-3\, D\vec{H}+\star(\ast_{g}d\vec{n}_{\vec{\Phi}}\wedge\vec{H})\rg]=\Im(q,\vec{h}_0)_{WP}
\ee
for some holomorphic quadratic differential $q$ associated to the fixed conformal class. $q$ plays here the role of a Lagrange multiplyer. Equation (\ref{I.14})
has been called {\it Constrained Willmore equation} (see \cite{BPP} for instance) but in order to avoid any ambiguity with the other constrained problems for the Willmore lagarngian (such as the Isoperimetric ratio for instance - see \cite{Sy}) we prefer to call equation (\ref{I.14}) the {\it Constrained-conformal Willmore} equation.

\medskip

Examples of solutions to (\ref{I.14}), which are not necessarily solutions\footnote{Surfaces of non-zero constant mean curvature in ${\R}^3$ which are Willmore have to be umbilic
and then coincide with a plane or a round sphere.}  to (\ref{I.12}) are given for instance by {\bf parallel mean curvature surfaces} : surfaces that
generalize to arbitrary codimensions the {\bf constant mean curvature equation} and that are characterized by the following condition
\be
\label{I.140a}
D\vec{H}=\pi_{\vec{n}}(d \vec{H})\equiv 0\quad.
\ee
Indeed, the Codazzi-Mainardi identity for a general conformal immersion $\vec{\Phi}$ of the disc $D^2$ reads (see \cite{Ri1})
\be
\label{I.140b}
e^{-2\la}\p_{\ov{z}}\lf(e^{2\la}\ {\vec{H}_0}\cdot\vec{H}\rg)=\vec{H}\cdot\p_{z}\vec{H}+{\vec{H}_0}\cdot\p_{\ov{z}}\vec{H}\quad,
\ee
where $z=x_1+ix_2$ and $\p_z:=2^{-1}(\p_{x_1}-i\p_{x_2})$. Since we are assuming (\ref{I.140a}) we have then
\be
\label{I.140c}
f(z):=e^{2\la}\ {\vec{H}_0}\cdot\vec{H}\quad\quad\mbox{ is holomorphic.}
\ee
In \cite{Ri1} it is proven that, for a general conformal immersion $\vec{\Phi}$ of the disc $D^2$ , one has
\be
\label{I.140d}
div\lf(\nabla\vec{H}-3\,\pi_{\vec{n}}(\nabla\vec{H})+\star(\nabla^\perp\vec{n}_{\vec{\Phi}}\wedge\vec{H})\rg)=-8\ \Re\lf(\p_{\ov{z}}\lf[\pi_{\vec{n}}(\p_z\vec{H})+ e^{2\la}\ {\vec{H}_0}\cdot\vec{H}\ e^{-2\la}
\p_{\ov{z}}\vec{\Phi}\rg]\rg)
\ee
Assuming (\ref{I.140a}), (\ref{I.140d}) becomes
\be
\label{I.140e}
div\lf(\nabla\vec{H}-3\,\pi_{\vec{n}}(\nabla\vec{H})+\star(\nabla^\perp\vec{n}_{\vec{\Phi}}\wedge\vec{H})\rg)=-8\ \Re\lf( f(z)\ \p_{\ov{z}}\lf[e^{-2\la}\ \p_{\ov{z}}\vec{\Phi}\rg]\rg)\quad.
\ee
For a general conformal immersion $\vec{\Phi}$ of the disc $D^2$ , one has (see \cite{Ri1})
\be
\label{I.140f}
\vec{H}_0=2\ \p_{z}\lf[e^{-2\la}\ \p_{z}\vec{\Phi}\rg]\quad.
\ee
Hence (\ref{I.140e}) becomes
\be
\label{I.140g}
div\lf(\nabla\vec{H}-3\,\pi_{\vec{n}}(\nabla\vec{H})+\star(\nabla^\perp\vec{n}_{\vec{\Phi}}\wedge\vec{H})\rg)=\Im\lf(4\,i\,f(z)\ \ov{\vec{H}_0}\rg)\quad, 
\ee
which is exactly the {\it constrained-conformal Willmore} equation (\ref{I.14}) written in complex coordinates.

\medskip

 If instead the critical point is a degenerate
point of the conformal constraint, it is proved in \cite{Ri3} that there exists a non trivial holomorphic quadratic differential $q$ such that
\be
\label{I.14a}
\Im(q,\vec{h}_0)_{WP}\equiv 0
\ee
in other words, $\vec{\Phi}$ is {\it isothermic}. 

We have proven in \cite{Ri3} (see propositions I.2 and I.3) that, if $\vec{\Phi}$ is 
isothermic , away from the zeros of $q$, 
there exists locally complex coordinates $z=x_1+ix_2$ in which the condition (\ref{I.8}) reads
\be
\label{I.14b}
\frac{\p}{\p x_1}\lf[e^{-2\la}\frac{\p\vec{\Phi}}{\p x_2}\rg]+\frac{\p}{\p x_2}\lf[e^{-2\la}\frac{\p\vec{\Phi}}{\p x_1}\rg]=0\quad.
\ee
where $e^\lambda=|\p_{x_1}\vec{\Phi}|=|\p_{x_2}\vec{\Phi}|$ is the conformal factor.

 Making a similar choice of conformal coordinates for the induced metric $g$ equation (\ref{I.14}) becomes.
\be
\label{I.15}
div\lf(\nabla\vec{H}-3\,\pi_{\vec{n}}(\nabla\vec{H})+\star(\nabla^\perp\vec{n}_{\vec{\Phi}}\wedge\vec{H})\rg)= Q\ \lf(\frac{\p}{\p x_1}\lf[e^{-2\la}\frac{\p\vec{\Phi}}{\p x_2}\rg]+\frac{\p}{\p x_2}\lf[e^{-2\la}\frac{\p\vec{\Phi}}{\p x_1}\rg]\rg)\quad.
\ee
where\footnote{By dilating these conformal coordinates
one can always make $Q=1$ in (\ref{I.15}) - except when $q=0$ of course - but we prefer to normalize the conformal coordinates for them not to degenerate as the Weil-Petersson norm of the Lagrange multiplier
  $|q|_{WP}$ would go either to $0$ or $+\infty$} $Q:=|q|_{WP}\in{\R}^+$ (The WP -norm is taken with respect to the constant scalar curvature metric of volume 1 on $\Sigma$). 

The {\bf Isothermic equation} (\ref{I.14b}) is an {\bf hyperbolic} equation whereas the {\bf Constrained-conformal Willmore equation} (\ref{I.15}) is an {\bf elliptic} one.
One passes from (\ref{I.15}) to (\ref{I.14b}) in particular when the norm of the Lagrange-multiplier goes to infinity $Q=\ep^{-2}\rightarrow +\infty$. Precisely in \cite{Ri3}
section IV we have proven the following result

\begin{Th}
\label{th-I.1}\cite{Ri3}
Let $\vec{\Phi}_k$ be a sequence of conformal immersion from $D^2$ into ${\R}^m$ satisfying asymptotically  the constrained-conformal equation :
\be
\label{I.16}
\begin{array}{l}
\ds div\lf(\nabla\vec{H}_k-3\,\pi_{\vec{n}_k}(\nabla\vec{H}_k)+\star(\nabla^\perp\vec{n}_{\vec{\Phi}_k}\wedge\vec{H}_k)\rg)- Q_k\ \lf(\frac{\p}{\p x_1}\lf[e^{-2\la_k}\frac{\p\vec{\Phi}_k}{\p x_2}\rg]+\frac{\p}{\p x_2}\lf[e^{-2\la_k}\frac{\p\vec{\Phi}_k}{\p x_1}\rg]\rg)\\[5mm]
\ds \quad\quad\longrightarrow 0\quad\quad \mbox{ strongly in }\quad(W^{2,2}\cap W^{1,\infty})^\ast
\end{array}
\ee
for some sequence $Q_k\in{\R}^+$. Assume 
\be
\label{I.17}
\|\la_k\|_{L^\infty(D^2)}+\|\nabla\vec{n}_{\vec{\Phi}_k}\|_{L^2(D^2)}\le C<+\infty\quad.
\ee
If
\[
\limsup_{k\rightarrow+\infty} Q_k<+\infty
\]
then, modulo extraction of a subsequence, $\vec{\Phi}_k$ converges weakly\footnote{in this case the weak $W^{2,2}_{loc}$ convergence should even be strong.} in $W^{2,2}_{loc}$ to a $C^\infty$ constrained conformal immersion (i.e. satisfying (\ref{I.15})
for some $Q\in {\R}^+$). 

\noindent Alternatively, if instead,
\[
\limsup_{k\rightarrow+\infty} Q_k=+\infty
\]
there exists a subsequence of $\vec{\Phi}_k$ converging weakly in $W^{2,2}_{loc}$ to a conformal lipschitz $W^{2,2}_{loc}$ isothermic immersion (i.e. satisfying (\ref{I.14b}))
\hfill $\Box$ 
\end{Th}
In this sense the {\bf isothermic surface equation} should be seen as an {\bf hyperbolic degeneracy} of the {\bf constrained conformal equation} which represents some
{\bf viscous approximation} of the first one. 
\begin{Rm}
\label{rm-I.2}
An interesting issue is to understand if the solution to (\ref{I.14b}) that are obtained as weak limits of the viscous approximation (\ref{I.15}) enjoy some additional 
regularity properties which are not shared with the arbitrary $W^{2,2}$ conformal solutions to   (\ref{I.14b}) .\hfill $\Box$
\end{Rm}
\subsection{Weak Global Isothermic Immersions.}

The previous result, theorem~\ref{th-I.1}, shows the importance of enlarging the class of {\it smooth global isothermic immersions} to a wider class
of {\it weak global isothermic immersions}. For analysis reasons it is also needed to enlarge the class of $C^1$ immersions while studying critical points to the Willmore functional
(\ref{I.10a}). In \cite{Ri3} the author introduced the framework of {\it weak immersion with finite total curvature} (or simply {\it weak immersions}).

Let $g_0$ be a reference smooth metric on $\Sigma$. One defines the Sobolev spaces $W^{k,p}(\Sigma,{\R}^m)$ of measurable maps from $\Sigma$ into 
${\R}^m$ in the following way
\[
W^{k,p}(\Sigma,{\R}^m)=\lf\{f\ \mbox{ meas. } {\Sigma}\rightarrow {\R}^m\mbox{ s.t. }\sum_{l=0}^k\int_{\Sigma}|\nabla^l f|_{g_0}^p\ dvol_{g_0}<+\infty\rg\}
\]
Since $\Sigma$ is assumed to be compact it is not difficult to see that this space is independent of the choice we have made of $g_0$.

\medskip

First we need to have a weak first fundamental form that is we need  $\vec{\Phi}^\ast g_{{\R}^m}$ to define an $L^\infty$ metric with a bounded inverse.
The last requirement is satisfied if we assume that $\vec{\Phi}$ is in $W^{1,\infty}(\Sigma)$ and if $d\vec{\Phi}$ has maximal rank 2 at every point with some uniform quantitative control
of ''how far'' $d\vec{\Phi}$ is from being degenerate : there exists $c_0>0$ s.t.
\be
\label{I.18}
|d\vec{\Phi}\wedge d\vec{\Phi}|_{g_0}\ge c_0>0 \quad.
\ee
where $d\vec{\Phi}\wedge d\vec{\Phi}$ is a 2-form on $\Sigma$ taking values into 2-vectors from ${\R}^m$ and given in local coordinates
by $2\,\p_x\vec{\Phi}\wedge\p_y\vec{\Phi}\ dx\wedge dy$.         
The condition (\ref{I.1}) is again independent of the choice of the metric $g_0$ .
For a Lipschitz immersion satisfying (\ref{I.1}) we can define the Gauss map as being the following measurable map in $L^\infty(\Sigma)$ taking values in the Grassmanian
of oriented $m-2$-planes in ${\R}^m$.
\[
\vec{n}_{\vec{\Phi}}:=\star\frac{\p_x\vec{\Phi}\wedge\p_y\vec{\Phi}}{|\p_x\vec{\Phi}\wedge\p_y\vec{\Phi}|}\quad.
\] 
We then introduce the space ${\mathcal E}_{\Sigma}$ of {\it weak immersions of $\Sigma$ 
with total finite curvature} as being the following space :
\[
\mathcal{E}_\Sigma:=\lf\{
\begin{array}{c}
\ds\vec{\Phi}\in W^{1,\infty}(\Sigma,{\R}^m)\quad\mbox{ s.t. } \vec{\Phi} \mbox{ satisfies }(\ref{I.18})\mbox{ for some }c_0\\[5mm]
\ds\mbox{ and }\quad\int_{\Sigma}|d\vec{n}|_g^2\ dvol_g<+\infty
\end{array}
\rg\}\quad .
\]
Where $g:=\vec{\Phi}^\ast g_{{\R}^m}$ is the pull back by $\vec{\Phi}$ of the flat  canonical metric $g_{{\R}^m}$ of ${\R}^m$ and $dvol_g$ is the volume form associated to $g$.

\medskip

The analysis of $\mathcal{E}_\Sigma$ shows that for completeness purposes (see \cite{Ri3}) one has to relax the fact that $\vec{\Phi}$ is globally an immersion
by requiring only that $\vec{\Phi}$ is an immersion away from finitely many points. We then define the space of {\it branched weak immersions with finite total curvature}
in the following way
\[
\mathcal{F}_\Sigma:=\lf\{
\begin{array}{c}
\ds\vec{\Phi}\in W^{1,\infty}(\Sigma,{\R}^m)\quad\mbox{ s.t. } \exists \ a_1\cdots a_N\in \Sigma \mbox{ s.t. } \\[5mm]
\forall\ K\mbox{ compact in }\Sigma\setminus\{a_1\cdot a_N\}\quad\vec{\Phi}\mbox{ satisfies (\ref{I.18})}\mbox{ on $K$ for some }c_0(K)>0\\[5mm]
\ds\mbox{ and }\quad\int_{\Sigma}|d\vec{n}|_g^2\ dvol_g<+\infty
\end{array}
\rg\}\quad .
\]

\medskip

It is proved in \cite{Ri3} (see also \cite{Ri1}) that any {\it weak immersion} $\vec{\Phi}$  in ${\mathcal E}_\Sigma$ defines a smooth conformal structure on $\Sigma$ : more precisely, following Toro, M\"uller-Sverak, H\'elein's works on immersions with finite total curvature one proves (see \cite{Ri1}) that  
for any $\vec{\Phi}\in {\mathcal E}_\Sigma$ and for any $p\in\Sigma$ there exists a neighborhood $U$ containing $p$ and a bilipshitz homeomorphism $\Psi$ from $D^2$
into $U$ such that
$\vec{\Phi}\circ\Psi$ satisfies the weak conformal condition 
\[
\lf\{
\begin{array}{l}
\ds \p_{x_1}(\vec{\Phi}\circ\Psi)\cdot  \p_{x_2}(\vec{\Phi}\circ\Psi)=0\quad\quad\mbox{ a. e. in }D^2\\[5mm]
\ds |\p_{x_1}(\vec{\Phi}\circ\Psi)|=|\p_{x_2}(\vec{\Phi}\circ\Psi)|\quad\quad\mbox{ a. e. in }D^2
\end{array}
\rg.
\]
moreover $\vec{\Phi}\circ\Psi$ is $W^{2,2}$ on $D^2$. Hence $\Sigma$ is equipped with a system of charts such that the  transition
functions satisfy the Cauchy-Riemann  conditions almost everywhere and thus are holomorphic. This defines the conformal structure induced by $\vec{\Phi}$. The same
can be done for any element of $\mathcal{F}_\Sigma$ using also Huber theorem about the conformal structure of a metric 
of finite total curvature on a closed surface minus finitely many points.

We can now give the definition of a weak global isothermic immersion as the natural extension of definition~\ref{df-I.1}.

\begin{Dfi}
\label{df-I.2}
Let $\Sigma^2$ be a closed two dimensional manifold. One says that a weak immersion $\vec{\Phi}$ in ${\mathcal E}_\Sigma$ (resp. a weak branched immersion in $\mathcal{F}_\Sigma$) is {\bf weakly global isothermic} if
there exists an holomorphic quadratic form $q$ of the riemann surface issued from $\Sigma^2$ equipped with the conformal structure defined by $\vec{\Phi}$ such that
\be
\label{I.19}
\Im(q,\vec{h}_0)_{WP}=0\quad.
\ee
where $\vec{h}_0$ is the Weingarten form of the immersion $\vec{\Phi}$ given by (\ref{I.9}).\hfill $\Box$
\end{Dfi}
\begin{Rm}
\label{rm-I.3}
Observe that for any $\vec{\Phi}$ in ${\mathcal E}_\Sigma$ the Weingarten $1-0\otimes 1-0$ form $h_0$ is a well defined $L^2$ section of
 $\wedge^{(1,0)}\Sigma\otimes\wedge^{(1,0)}\Sigma$ and therefore the function $\Im(q,\vec{h}_0)_{WP}$ is a well efined $L^2$ function on $\Sigma$
 for any holomorphic quadratic form $q$.
 \end{Rm}
 The following characterization of weak global isothermic immersion has been given in \cite{Ri3} (proposition I.3).
 \begin{Prop}
 \label{pr-I.3}
 A weak immersion $\vec{\Phi}$ is global isothermic if and only if around every point there exists a $L^2$ ${\R}^m$ valued map $\vec{L}$ such that
 the following two conditions are satisfied
 \be
 \label{I.19a}
 \lf\{
 \begin{array}{l}
\ds d\vec{\Phi}\cdot d\vec{L}:=[\p_{x_1}\vec{\Phi}\cdot\p_{x_2}\vec{L}-\p_{x_2}\vec{\Phi}\cdot\p_{x_1}\vec{L}]\ dx_1\wedge dx_2=0\\[5mm]
\ds d\vec{\Phi}\wedge d\vec{L}:=[\p_{x_1}\vec{\Phi}\wedge\p_{x_2}\vec{L}-\p_{x_2}\vec{\Phi}\wedge\p_{x_1}\vec{L}]\ dx_1\wedge dx_2=0
\end{array}
\rg.
 \ee
$\vec{L}$ is called a {\bf Darboux transform} of $\vec{\Phi}$. \hfill$\Box$
 \end{Prop} 
 An elementary observation shows that property (\ref{I.19a}) is invariant under the action of transformations that preserves angles infinitesimally in ${\R}^m$. 
 From this observation we deduce the following fundamental property.
 \begin{Prop}
 \label{pr-I.4}
 Let $\vec{\Phi}$ be a weak isothermic immersion of $\mathcal{E}_\Sigma$ (resp. weak branched isothermic immersion of $\mathcal{F}_\Sigma$). Let $\Xi$ be a conformal
 transformation of ${\R}^m\cup\{\infty\}$. Then $\vec{\Xi}\circ\vec{\Phi}$ is still a weak isothermic immersion of $\mathcal{E}_\Sigma$ (resp. weak branched isothermic immersion of $\mathcal{F}_\Sigma$).\hfill $\Box$
 \end{Prop}
 In \cite{BR} the following proposition is proved
 \begin{Prop}
 \label{pr-I.5}
 A weak immersion $\vec{\Phi}\in{\mathcal E}_\Sigma$ is constrained-conformal Willmore if and only if, around every point, there exists an $L^2$ ${\R}^m$-valued map $\vec{L}$ such that
 the following two conditions are satisfied
 \be
 \label{I.19b}
 \lf\{
 \begin{array}{l}
\ds d\vec{\Phi}\cdot d\vec{L}:=[\p_{x_1}\vec{\Phi}\cdot\p_{x_2}\vec{L}-\p_{x_2}\vec{\Phi}\cdot\p_{x_1}\vec{L}]\ dx_1\wedge dx_2=0\\[5mm]
\ds d\vec{\Phi}\wedge d\vec{L}:=[\p_{x_1}\vec{\Phi}\wedge\p_{x_2}\vec{L}-\p_{x_2}\vec{\Phi}\wedge\p_{x_1}\vec{L}]\ dx_1\wedge dx_2=2\ (-1)^m\ d(\star(\vec{n}\res\vec{H}))\wedge\res d\vec{\Phi}\quad,
\end{array}
\rg.
 \ee
 where  $\res$ is the standard contraction operator between  a $p-$vectors and a $q-$vectors  $(p\ge q)$ given by
\[
\forall\,\vec{a}\in\wedge^p{\R}^m\ ,\quad\forall\,\vec{b}\in\wedge^q{\R}^m\, ,\quad\forall\,\vec{c}\in\wedge^{p-q}{\R}^m\quad\quad<\vec{a}\res\vec{b},\vec{c}>=<\vec{a},\vec{b}\wedge\vec{c}>\quad.
\] 
and 
\[
d(\star(\vec{n}\res\vec{H}))\wedge\res d\vec{\Phi}:=[\p_{x_1}(\star(\vec{n}\res\vec{H}))\res \p_{x_2}\vec{\Phi}-\p_{x_2}(\star(\vec{n}\res\vec{H}))\res \p_{x_1}\vec{\Phi}]
\ dx_1\wedge dx_2\quad.
\]
\hfill $\Box$
 \end{Prop}
In \cite{Ri2} it is proven that  weak immersion $\vec{\Phi}\in{\mathcal E}_\Sigma$ which are {\it constrained-conformal Willmore} are in fact $C^\infty$.

\medskip

{\bf Minimal surfaces} in ${\R}^m$ -satisfying $\vec{H}=0$ - clearly solve (\ref{I.13}). This means that they are Willmore and, a fortiori, they are special cases of {\it constrained-conformal Willmore}. Therefore, from proposition~\ref{pr-I.5},
they satisfy (\ref{I.19b}). But since $\vec{H}=0$ the right hand side of (\ref{I.19b}) is zero. Thus minimal surfaces are also satisfying (\ref{I.19a}) and are then {\bf isothermic}.

\medskip

More generally {\bf parallel mean curvature surfaces}, surfaces satisfying (\ref{I.140a}), are also {\bf constrained-conformal Willmore} and not necessarily Willmore, as we proved
in the previous subsection, and they are also {\bf isothermic}.
Indeed, it is proven in \cite{BR} (equation ({II.6}) )that, in conformal coordinates,
\be
\label{I.190c}
\p_{x_1}(\star(\vec{n}\res\vec{H}))\res \p_{x_2}\vec{\Phi}-\p_{x_2}(\star(\vec{n}\res\vec{H}))\res \p_{x_1}\vec{\Phi}=(-1)^{m-1}\ \nabla\vec\Phi\wedge\nabla\vec{H}
\ee
For {\it parallel mean curvature surfaces}, which satisfy (\ref{I.140a}), we have
\be
\label{I.190d}
\nabla\vec\Phi\wedge\nabla\vec{H}=\nabla\vec\Phi\wedge\pi_T(\nabla\vec{H})=\nabla^\perp\vec{H}\cdot\nabla\vec{\Phi}\ \vec{e}_1\wedge\vec{e}_2=-2\,
div(\vec{H}\cdot\nabla^\perp\vec{\Phi})\ \vec{e}_1\wedge\vec{e}_2=0\quad.
\ee

\medskip

Other examples of weak isothermic immersions which are not smooth and then not necessarily {\it constrained-conformal Willmore} are easy to produce : take a non necessarily smooth simple closed lipshitz curve $\gamma\ :\ S^1\rightarrow {\R}^2$ such that
 \[
 \int_{S^1}\kappa^2\ dl<+\infty
 \]
 where $\kappa$ is the curvature distribution of that curve and $dl$ the length 1-form on $S^1$ induced by the immersion $\gamma$. Identify the plane ${\R}^2$ with the vertical
 plane in ${\R}^3$ given by $\{x_2=0\}$ and rotate that curve around the $x_3$ vertical axis. One proves that this generates a {\it weak global isothermic immersion} : {\bf axially symmetric surfaces} are {\bf isothermic}.
 We saw in proposition~\ref{pr-I.4} that being isothermic is a conformally invariant property  and therefore any composition of the obtained axially surface with a diffeomorphism of ${\R}^3$ generates another isothermic surface.

It is proven in \cite{Ri3} (see the proof of lemma III.1) that the space of weak immersion ${\mathcal E}_\Sigma$ of controlled conformal class has a nice weak closure property modulo
renormalization and branched points. Precisely one has the following weak closure lemma.
\begin{Lm}
\label{lm-I.1} \cite{Ri3}
Let $\Sigma$ be a closed two-dimensional manifold. Let $\vec{\Phi}_k$ be a sequence of elements in ${\mathcal E}_\Sigma$ such that $W(\Phi_k)$ is uniformly bounded.
Assume that the conformal class of the conformal structure $c_k$ (i.e. complex structure of $\Sigma$) defined by $\vec{\Phi}_k$ remains in a compact subspace of the Moduli space of $\Sigma$. Then, modulo
extraction of a subsequence, the sequence $c_k$ converges to a smooth limiting complex structure $c_\infty$ ; and there
exist a sequence of Lipschitz diffeomorphisms $f_k$ of $\Sigma$  such that $\vec{\Phi}_k\circ f_k$ is conformal from $(\Sigma,c_k)$ into ${\R}^m$. Moreover, there exists a sequence $\Xi_k$ 
of conformal diffeomorphisms of ${\R}^m\cup\{\infty\}$ and at most finitely many points $\{a_1,\ldots, a_N\}$ such that  
\be
\label{II.0}
\limsup_{k\rightarrow +\infty} {\mathcal H}(\Xi_k\circ\vec{\Phi}_k\circ f_k(\Sigma))<+\infty\quad\quad, \quad\quad\Xi_k\circ\vec{\Phi}_k\circ f_k(\Sigma)\subset B_R(0)
\ee
for some $R>0$ independent of $k$, and
\be
\label{II.1}
\vec{\xi}_k:=\Xi_k\circ\vec{\Phi}_k\circ f_k\;\rightharpoonup \;\vec{\xi}_\infty\quad\quad\mbox{ weakly in }\:(W^{2,2}_{loc}\cap W^{1,\infty}_{loc})^\ast(\Sigma\setminus\{a_1,\ldots, a_N\})\quad.
\ee
The convergences are understood with respect to $h_k$, which is the constant scalar curvature metric of unit volume attached to the conformal structure $c_k$. \\
Furthermore, there holds
\be
\label{II.2}
\forall\: K\mbox{ compact subset of }\Sigma\setminus\{a_1,\ldots, a_N\}\quad\limsup_{k\rightarrow+\infty}\|\log|d\vec{\xi}_k|_{h_k}\|_{L^\infty(K)}<+\infty\quad.
\ee
Finally, $\vec{\xi}_\infty$ is an element of ${\mathcal F}_\Sigma$, a weak immersion of $\Sigma\setminus\{a_1,\ldots, a_N\}$,
and conformal from $(\Sigma,c_\infty)$ into ${\R}^m$.\hfill$\Box$
\end{Lm}
Following the arguments of \cite{Ri3} proof of lemma IV.1 one establishes the following weak closure result for weak isothermic immersions.
\begin{Th}
\label{th-I.2} \cite{Ri3}
Let $\Sigma$ be a closed two-dimensional manifold. Let $\vec{\Phi}_k$ be a sequence of {\bf weak global isothermic immersions} such that $W(\Phi_k)$ is uniformly bounded.
Assume that the conformal classes $c_k$ defined by $\vec{\Phi}_k$ converge to a limiting structure $c_\infty$ in the Moduli space of $\Sigma$. Then, modulo extraction of a subsequence, there exists a sequence of Lipschitz diffeomorphisms $f_k$ of $\Sigma$ and a sequence $\Xi_k$ 
of conformal diffeomorphisms of ${\R}^m\cup\{\infty\}$ such that $\vec{\xi}_k:=\Xi_k\circ\vec{\Phi}_k\circ f_k$ is a weak conformal isothermic immersion converging weakly
in $W^{2,2}_{loc}$ on $\Sigma$ minus finitely many points to $\vec{\xi}_\infty$ a , possibly branched at these points, conformal {\bf weak global isothermic immersion} for the limiting conformal structure $c_\infty$ on $\Sigma$.\hfill $\Box$
\end{Th}
\subsection{Weakly converging smooth global isothermic immersions - Main result.}

The goal of the present paper is to present a result regarding the lack of strong compactness and the geometric structure of the defect measure for sequences of smooth global isothermic immersions weakly converging to another smooth global isothermic immersion. Our main result is the following

\begin{Th}
\label{th-I.3}
Let $\Sigma$ be a closed two-dimensional manifold. Let $\vec{\Phi}_k$ be a sequence of smooth {\bf global isothermic immersions} such that $W(\Phi_k)$ is uniformly bounded.
Assume that the conformal classes $c_k$ defined by $\vec{\Phi}_k$ converge to a limiting structure $c_\infty$ in the Moduli space of $\Sigma$. Then, modulo extraction of a subsequence, there exists a sequence of Lipschitz diffeomorphisms $f_k$ of $\Sigma$ and a sequence $\Xi_k$ 
of conformal diffeomorphisms of ${\R}^m\cup\{\infty\}$ and finitely many points $\{a^1\cdots a^n\}$ such that $\vec{\xi}_k:=\Xi_k\circ\vec{\Phi}_k\circ f_k$ is a conformal global isothermic immersion satisfying
\be
\label{I.20a}
\vec{\xi}_k\rightharpoonup\vec{\xi}_\infty\quad\quad\mbox{ weakly in}\quad W^{2,2}_{loc}(\Sigma\setminus\{a^1\cdots a^n\})
\ee
where $\vec{\xi}_\infty$ a weak, possibly branched at the $a^i$, conformal {\bf weak global isothermic immersion} for the limiting conformal structure $c_\infty$ on $\Sigma$. If moreover $\vec{\xi}_\infty$ is {\bf smooth} away from the points $a^i$ then the following convergence holds
\be
\label{I.20}
|d\vec{n}_{\vec{\xi}_k}|^2_{g_k}\ dvol_{g_k}\rightharpoonup|d\vec{n}_{\vec{\xi}_\infty}|^2_{g_\infty}\ dvol_{g_\infty}+\nu+\sum_{i=1}^n\al^i\ \delta_{a^i}\ dvol_{g_\infty}\quad\quad\mbox{ weakly in }\quad{\mathcal M}(\Sigma)
\ee
where ${\mathcal M}(\Sigma)$ is the space of Radon measures on $\Sigma$ and $\nu$, the non atomic part of the defect measure, satisfies the following condition : around every point
different from the $a^i$ there exists a conformal coordinate chart $z=x_1+ix_2$ such that, simultaneously the following holds 
\be
\label{I.21}
0=\Im({\vec{H}_0})=-\frac{1}{2}\pi_{\vec{n}}\lf(\frac{\p}{\p x_1}\lf[e^{-2\la}\frac{\p\vec{\xi}_\infty}{\p x_2}\rg]+\frac{\p}{\p x_2}\lf[e^{-2\la}\frac{\p\vec{\xi}_\infty}{\p x_1}\rg]\rg)
\ee 
where $\vec{H}_0$ is the expression in the $z$ coordinates of the Weingarten form $\vec{h}_0=\vec{H}_0\ dz\otimes dz$, and
\be
\label{I.22}
\nu=\nu_1(x_1)\wedge dx_2+dx_1\wedge\nu_2(x_2)
\ee
where $\nu_i(x_i)$ are Radon measures on the $x_i$ axis and $\nu_i(x_i)\wedge dx_{i+1}$ is the product of this Radon measure with the Lebesgue measure on the $x_{i+1}$ axis.
\hfill $\Box$
\end{Th}
\begin{Rm}
\label{rm-I.4}
In codimension 1 the coordinates directions in which (\ref{I.21}) happens are {\bf principal directions}. The theorem says that the defect measure associated to the lack of strong compactness
of the sequence of isothermic immersions  {\bf ''propagates'' uniformly along principal directions}, modulo possible concentration points. \hfill $\Box$
\end{Rm}
\begin{Rm}
\label{rm-I.5}
The result is optimal in the sense that it is not difficult to
produce examples where (\ref{I.22}) indeed happens. Consider a family of simple closed curves in the plane of fixed length, such that, the normal parametrization, $\gamma_k(s)$, weakly converges in $W^{2,2}(S^1)$ with a non zero defect measure $\mu(s)$
\[
|\ddot{\gamma}_k|^2(s)\ ds\rightharpoonup |\ddot{\gamma}_\infty|^2(s)\ ds+\mu(s)
\]
By identifying the 2-plane with the vertical plane in ${\R}^3$ given by $\{x_2=0\}$ and by rotating the sequence of curves around the $x_3$ axis we obtain a weakly converging family
of isothermic surfaces with a non zero defect measure satisfying (\ref{I.22}).\hfill $\Box$
\end{Rm}

\section{Entropies for Isothermic Surfaces.}
\reset
One of the main tool for proving theorem~\ref{th-I.3} is the computation of entropies for isothermic surfaces. Precisely the goal of the present section is to establish the following
proposition.
\begin{Prop}
\label{pr-II.1}
Let $\vec{\Phi}$ be a smooth conformal immersion of $D^2$ into ${\R}^m$ satisfying 
\be
\label{II.4}
\frac{\p}{\p x_1}\lf[e^{-2\la}\frac{\p\vec{\Phi}}{\p x_2}\rg]+\frac{\p}{\p x_2}\lf[e^{-2\la}\frac{\p\vec{\Phi}}{\p x_1}\rg]=0\quad.
\ee
where $e^\lambda=|\p_{x_1}\vec{\Phi}|=|\p_{x_2}\vec{\Phi}|$ is the conformal factor. Then the following conservation laws hold
\be
\label{II.5}
\lf\{
\begin{array}{l}
\ds\frac{\p}{\p x_1}\lf[ \lf(\frac{\p\vec{n}_\Phi}{\p x_2}\res\vec{e}_2\rg)^2+\lf|\frac{\p\la}{\p{x_1}}\rg|^2-\lf|\frac{\p\la}{\p{x_2}}\rg|^2\rg]+\frac{\p}{\p {x_2}}\lf[2\ \frac{\p\la}{\p{x_1}}\ \frac{\p\la}{\p{x_2}}\rg]=0\\[5mm]
\ds\frac{\p}{\p x_2}\lf[ \lf(\frac{\p\vec{n}_\Phi}{\p x_1}\res\vec{e}_1\rg)^2+\lf|\frac{\p\la}{\p{x_2}}\rg|^2-\lf|\frac{\p\la}{\p{x_1}}\rg|^2\rg]+\frac{\p}{\p {x_1}}\lf[2\ \frac{\p\la}{\p{x_1}}\ \frac{\p\la}{\p{x_2}}\rg]=0
\end{array}
\rg.
\ee
where $\vec{e}_i$ is the unit orthonormal Coulomb frame of $\vec{\Phi}_\ast T D^2$ given by $\vec{e}_i:=e^{-\la}\,\p_{x_i}\vec{\Phi}$ and
 $\res$ is the following standard contraction operator between  a $p-$vector and a $q-$vector  $(p\ge q)$ giving a $p-q-$vector
\[
\forall\,\vec{a}\in\wedge^p{\R}^m\ ,\quad\forall\,\vec{b}\in\wedge^q{\R}^m\, ,\quad\forall\,\vec{c}\in\wedge^{p-q}{\R}^m\quad\quad<\vec{a}\res\vec{b},\vec{c}>=<\vec{a},\vec{b}\wedge\vec{c}>\quad.
\] 

\hfill $\Box$
\end{Prop}
\begin{Rm}
\label{rm-II.1}
The proof of  proposition~\ref{pr-II.1} we give below is using the smoothness of the isothermic immersion and, a-priori (\ref{II.5}) does not necessarily hold for general
isothermic weak immersion in $\mathcal{E}_\Sigma$.\hfill $\Box$
\end{Rm}
\noindent{\bf Proof of proposition~\ref{pr-II.1}.}

A classical computation (see \cite{BR}) gives
\be
\label{II.6}
\vec{H}_0=2\ \p_z\lf[e^{-\la}\ \vec{e}_z\rg]
\ee
where $\p_z:=2^{-1}(\p_{x_1}-i\p_{x_2})$ and $\vec{e}_z:=2^{-1}(\vec{e}_1-i\vec{e}_2)$.  Observe that this identity implies
\be
\label{II.7}
\Im(\vec{H}_0)= 2^{-1} \frac{\p}{\p x_1}\lf[e^{-2\la}\frac{\p\vec{\Phi}}{\p x_2}\rg]+2^{-1}\frac{\p}{\p x_2}\lf[e^{-2\la}\frac{\p\vec{\Phi}}{\p x_1}\rg]\quad.
\ee
Our assumption is then equivalent to 
\be
\label{II.8}
\Im(\vec{H}_0)= 0\quad.
\ee
Since $\pi_{\vec{n}}(\vec{H}_0)=\vec{H}_0$, where $\pi_{\vec{n}}$ denotes the orthogonal projection onto the $m-2$ plane perpendicular to $\vec{e}_1$ and $\vec{e}_2$,
we deduce from (\ref{II.7}) and (\ref{II.8}) that
\be
\label{II.9}
0=\pi_{\vec{n}}\lf(\frac{\p}{\p x_1}\lf[e^{-2\la}\frac{\p\vec{\Phi}}{\p x_2}\rg]+\frac{\p}{\p x_2}\lf[e^{-2\la}\frac{\p\vec{\Phi}}{\p x_1}\rg]\rg)\quad,
\ee
which itself implies
\be
\label{II.10}
\pi_{\vec{n}}\lf(\frac{\p^2\vec{\Phi}}{\p x_1\p x_2}\rg)=0\quad.
\ee
Observe that we have
\[
\p_{x_2}\vec{n}\res\p_{x_1}\vec{\Phi}=\p_{x_2}(\vec{n}\res\p_{x_1}\vec{\Phi})-\vec{n}\res\p^2_{x_1x_2}\vec{\Phi}=-\vec{n}\res\p^2_{x_1x_2}\vec{\Phi}=\p_{x_1}(\vec{n}\res\p_{x_2}\vec{\Phi})-\vec{n}\res\p^2_{x_1x_2}\vec{\Phi}=\p_{x_1}\vec{n}\res\p_{x_2}\vec{\Phi}\quad.
\]
where we have used that $\vec{n}\res\p_{x_1}\vec{\Phi}=0$ and $\vec{n}\res\p_{x_2}\vec{\Phi}$. Inserting (\ref{II.10}) in this identity gives
\be
\label{II.11}
\frac{\p\vec{n}}{\p x_2}\res\frac{\p\vec{\Phi}}{\p x_1}=\frac{\p\vec{n}}{\p x_1}\res\frac{\p\vec{\Phi}}{\p x_2}=0\quad.
\ee
We have
\be
\label{II.12}
\begin{array}{l}
\ds\p_{x_1}(\p_{x_2}\vec{n}\res\vec{e}_2)=\p_{x_1}(e^{-\la}\ \p_{x_2}\vec{n}\res\p_{x_2}\vec{\Phi})\\[5mm]
\ds=-\p_{x_1}\la\ \p_{x_2}\vec{n}\res\vec{e}_2+e^{-\la}\ \p^2_{x_1 x_2}\vec{n}\res\p_{x_2}\vec{\Phi}+e^{-\la}\ \p_{x_2}\vec{n}\res\p^2_{x_1x_2}\vec{\Phi}\\[5mm]
\ds=-\p_{x_1}\la\ \p_{x_2}\vec{n}\res\vec{e}_2-e^{-\la}\ \p_{x_1}\vec{n}\res\p^2_{x_2^2}\vec{\Phi}+e^{-\la}\ \p_{x_2}\vec{n}\res\pi_T(\p^2_{x_1x_2}\vec{\Phi})
\end{array}
\ee
where we have used (\ref{II.10}) and (\ref{II.11}). In one hand we have 
\be
\label{II.13}
\begin{array}{l}
\ds\pi_T(\p^2_{x_1x_2}\vec{\Phi})=2^{-1}\,e^{-\la}\ \lf[\p_{x_2}(|\p_{x_1}\vec{\Phi}|^2)\ \vec{e}_1+\p_{x_1}|\p_{x_2}\vec{\Phi}|^2\ \vec{e}_2\rg]\\[5mm]
\ds\quad\quad=\p_{x_2}\la\ \p_{x_1}\vec{\Phi}+\p_{x_1}\la\ \p_{x_2}\vec{\Phi}
\end{array}
\ee
Thus using (\ref{II.11}) we have
\be
\label{II.14}
e^{-\la}\ \p_{x_2}\vec{n}\res\pi_T(\p^2_{x_1x_2}\vec{\Phi})=\p_{x_1}\la\ e^{-\la}\ \p_{x_2}\vec{n}\res\p_{x_2}\vec{\Phi}=\p_{x_1}\la\ \p_{x_2}\vec{n}\res\vec{e}_2\quad,
\ee
and (\ref{II.12}) becomes
\be
\label{II.15}
\ds\p_{x_1}(\p_{x_2}\vec{n}\res\vec{e}_2)=-e^{-\la}\ \p_{x_1}\vec{n}\res\p^2_{x_2^2}\vec{\Phi}
\ee
In the other hand
\be
\label{II.16}
<\p_{x_1}\vec{n}\res\pi_{\vec{n}}(\p^2_{x_2^2}\vec{\Phi}), \p_{x_2}\vec{n}\res\vec{e}_2>=<\p_{x_1}\vec{n},\pi_{\vec{n}}(\p^2_{x_2^2}\vec{\Phi})\wedge(\p_{x_2}\vec{n}\res\vec{e}_2)>=0
\ee
Indeed, if $m=3$ $\p_{x_1}\vec{n}$ is perpendicular to the vector $\vec{n}$ to which $\pi_{\vec{n}}(\p^2_{x_2^2}\vec{\Phi})$ is parallel and, in the case when $m>3$, one easily verifies that
\[
(\p_{x_2}\vec{n}\res\vec{e}_2)\res\vec{e}_i=0\quad\quad\mbox{ for }i=1,2\quad,
\]
thus $\pi_{\vec{n}}(\p^2_{x_2^2}\vec{\Phi})\wedge(\p_{x_2}\vec{n}\res\vec{e}_2)$ is paralel to $\vec{n}$ which proves (\ref{II.16}).

\medskip

Combining now (\ref{II.15}) and (\ref{II.16}) we obtain
\[
<\p_{x_1}(\p_{x_2}\vec{n}\res\vec{e}_2),\p_{x_2}\vec{n}\res\vec{e}_2>=-e^{-\la}\ <\p_{x_1}\vec{n}\res\pi_T(\p^2_{x_2^2}\vec{\Phi}),\p_{x_2}\vec{n}\res\vec{e}_2>\quad.
\]
Using two more times (\ref{II.11}) this gives
\be
\label{II.17}
\begin{array}{l}
\ds<\p_{x_1}(\p_{x_2}\vec{n}\res\vec{e}_2),\p_{x_2}\vec{n}\res\vec{e}_2>=-e^{-\la}\ <\p_{x_1}\vec{n}\res\vec{e}_1,\p_{x_2}\vec{n}\res\vec{e}_2>\ <\vec{e}_1,\p^2_{x_2^2}\vec{\Phi}>\\[5mm]
\ds\quad\quad\quad=\lf[<\p_{x_1}\vec{n}\res\vec{e}_1,\p_{x_2}\vec{n}\res\vec{e}_2>-|\p_{x_1}\vec{n}\res\vec{e}_2|^2\rg]\ \p_{x_1}\la\\[5mm]
\ds\quad\quad\quad= e^{2\la}\ K\ \p_{x_1}\la=-\Delta\la\ \p_{x_1}\la=-\p_{x_1}(|\p_{x_1}\la|^2/2)-\p_{x_2}(\p_{x_1}\la\,\p_{x_2}\la)+\p_{x_1}(|\p_{x_2}\la|^2/2)
\end{array}
\ee
where $K$ is the Gauss curvature and where we have used the Liouville equation. (\ref{II.17}) gives the first equation of (\ref{II.5}). The second equation is established
in a similar way. The proof of proposition~\ref{pr-II.1} is complete.\hfill $\Box$

\section{A lemma in Compensation Compactness Theory}
\reset
In order to prove the main theorem~\ref{th-I.3} we shall need a compactness result related to some quantites present in the
expressions (\ref{II.5})  of the entropies. This result is based
on a compensation phenomenon observed first in \cite{De} (see also \cite{Ge} and \cite{EM}) in the framework of the analysis
of 2-dimensional perfect incompressible fluids. 

\begin{Lm}
\label{lm-III.1}
Let $\al_k$ and $\beta_k$ be two sequences of functions in $W^{1,2}(D^2,{\R})$
\be
\label{III.11}
\limsup_{k\rightarrow +\infty}\|\nabla\al_k\|_{L^2(D^2)}+\|\nabla\beta_k\|_{L^2(D^2)}<+\infty
\ee
Let $\varphi_k$ be the sequence of solutions in $W^{1,2}(D^2,{\R})$ of
\be
\label{III.12}
\lf\{
\begin{array}{l}
\ds\Delta\varphi_k=\p_{x_1}\al_k\,\p_{x_2}\beta_k-\p_{x_2}\al_k\,\p_{x_1}\beta_k\quad\quad\mbox{ in }D^2\\[5mm]
\ds\varphi_k=0\quad\quad\quad\mbox{ on }\partial D^2
\end{array}
\rg.
\ee
Then there exists a subsequence $\varphi_{k'}$ and two Radon measures $\mu$ and $\nu$ such that
\be
\label{III.13}
\lf\{
\begin{array}{l}
|\p_{x_1}\varphi_{k'}|^2-|\p_{x_2}\varphi_{k'}|^2\rightharpoonup|\p_{x_1}\varphi_{\infty}|^2-|\p_{x_2}\varphi_{\infty}|^2+\mu\quad\quad\mbox{ in }\quad{\mathcal D}'(D^2)\\[5mm]
\p_{x_1}\varphi_{k'}\ \p_{x_2}\varphi_{k'}\rightharpoonup\p_{x_1}\varphi_{\infty}\ \p_{x_2}\varphi_{\infty}+\nu\quad\quad\mbox{ in }\quad{\mathcal D}'(D^2)
\end{array}
\rg.
\ee
where
\be
\label{III.14}
\lf\{
\begin{array}{l}
\ds\Delta\varphi_\infty=\p_{x_1}\al_\infty\,\p_{x_2}\beta_\infty-\p_{x_2}\al_\infty\,\p_{x_1}\beta_\infty\quad\quad\mbox{ in }\quad D^2\\[5mm]
\ds\varphi_\infty=0\quad\quad\quad\mbox{ on }\quad\partial D^2
\end{array}
\rg.
\ee
and $\al_\infty$ (resp. $\beta_\infty$) is the weak limit in $W^{1,2}$ of $\al_{k'}$ (resp. $\beta_{k'}$). Moreover both $\mu$ and $\nu$ are atomic inside $D^2$ :
there exists $p_i\in D^2$ for $i\in {\N}$, and $q_j\in D^2$ for $j\in {\N}$ such that
\be
\label{III.15}
\mu=\sum_{i\in {\N}}c_i\, \delta_{p_i}\quad\quad\mbox{ and }\quad\quad\nu=\sum_{j\in {\N}}d_j\, \delta_{q_j}\quad\quad\mbox{in }{\mathcal D}'(D^2)\quad,
\ee
where
\be
\label{III.16}
\sum_{i\in {\N}}|c_i|=|\mu|(D^2)<+\infty\quad\quad\mbox{ and }\quad\quad\sum_{j\in {\N}}|d_j|=|\nu|(D^2)<+\infty\quad.
\ee
\hfill $\Box$
\end{Lm}
\noindent{Proof of lemma~\ref{lm-III.1}.}

Let $\ti{\al}_k$ and $\ti{\beta}_k$ be two Whitney type extension on the whole plane ${\R}^2$ of respectively $\al_k$ and $\beta_k$ satisfying
\be
\label{III.17}
\|\nabla\ti{\al}_k\|_{L^2({\R}^2)}\le C\ \|\nabla\al_k\|_{L^2(D^2)}\quad\quad\mbox{ and }\quad\quad\|\nabla\ti{\beta}_k\|_{L^2({\R}^2)}\le C\ \|\nabla\beta_k\|_{L^2(D^2)}
\ee
where $C$ is independent of the two sequences $\al_k$ and $\beta_k$ (take for instance in ${\R}^2\setminus D^2$ respectively $\ti{\al}_k(x):=\al(x/|x|^2)$ and $\ti{\beta}_k(x):=\beta_k(x/|x|^2)$.
Introduce
\be
\label{III.17a}
\ti{\varphi}_k:=\frac{1}{2\pi}\ \log|x|\ \ast\ \lf[\p_{x_1}\ti{\al}_k\,\p_{x_2}\ti{\beta}_k-\p_{x_2}\ti{\al}_k\,\p_{x_1}\ti{\beta}_k\rg]\quad.
\ee

\medskip

From Wente theorem (see \cite{We} and the exposition in \cite{He}) we know that both $\varphi_k$ and $\ti{\varphi}_k$ are uniformly bounded in $W^{1,2}\cap L^\infty$
and we have in particular
\be
\label{III.18}
\|\varphi_k\|_{L^\infty({\R}^2)}+\|\ti{\varphi}_k\|_{L^\infty({\R}^2)}+\|\nabla\varphi_k\|_{L^2({\R}^2)}+\|\nabla\ti{\varphi}_k\|_{L^2({\R}^2)}\le C\  \|\nabla\al_k\|_{L^2(D^2)}\ \|\nabla\beta_k\|_{L^2(D^2)}
\ee
Hence the difference $v_k=\ti{\varphi}_k-\varphi_k$, which is harmonic in $D^2$, is strongly precompact in every $C^l_{loc}(D^2)$ for $l\in{\N}$ and since we don't care about
concentration of the measures at the boundary $\p D^2$, it suffices to prove the results of the lemma (identities (\ref{III.13}...\ref{III.16}) for $\ti{\varphi}_k$ in $D^2$, this will imply
the corresponding identities for $\varphi_k$ in $D^2$

\medskip

We present the proof of the lemma for the quantity $\p_{x_1}\ti{\varphi}_k\,\p_{x_2}\ti{\varphi}_k$ (the proof for the other quantity $|\p_{x_1}\ti{\varphi}_k|^2-|\p_{x_2}\ti{\varphi}_k|^2$
being identical).

\medskip

To shorten a bit the notation we write
\[
\om_k:=\p_{x_1}\ti{\al}_k\,\p_{x_2}\ti{\beta}_k-\p_{x_2}\ti{\al}_k\,\p_{x_1}\ti{\beta}_k=\Delta\ti{\varphi}_k\quad.
\]
Because of the uniform bounds given by (\ref{III.18}) combined with the assumption (\ref{III.11}),
we can extract a subsequence still denoted $\ti{\varphi}_k$ such that
\[
\ti{\varphi}_k\rightharpoonup\ti{\varphi}_\infty\quad\quad\mbox{ weakly in }\quad W^{1,2}({\R}^2)
\]
and, due to the jacobian structure, we can pass to the limit in (\ref{III.17a}) :
\[
\ti{\varphi}_\infty:=\frac{1}{2\pi}\ \log|x|\ \ast\ \lf[\p_{x_1}\ti{\al}_\infty\,\p_{x_2}\ti{\beta}_\infty-\p_{x_2}\ti{\al}_\infty\,\p_{x_1}\ti{\beta}_\infty\rg]
\]
where $\ti{\al}_\infty$ and $\ti{\beta}_\infty$ are weak $W^{1,2}$-limits of respectively $\ti{\al}_k$ and $\ti{\beta}_k$. Moreover we can also ensure that
\[
\p_{x_1}\ti{\varphi}_k\, \p_{x_2}\ti{\varphi}_k\rightharpoonup \gamma\quad\quad\mbox{ weakly  in }\quad{\mathcal M}(D^2)
\]
where ${\mathcal M}(D^2)$ denotes the space of Radon measures. It remains now to identify the Radon measure $\gamma$.

\medskip

Let $\psi$ be an arbitrary function in $C^\infty_0(D^2)$, denoting by $\Delta^{-1}$ the convolution with $(2\pi)^{-1}\ \log\,|x|$ we have
\be
\label{III.19}
\begin{array}{l}
\ds\int_{D^2}\psi(x)\ \p_{x_1}\ti{\varphi}_k\, \p_{x_2}\ti{\varphi}_k\ dx=\int_{{\R}^2}\psi(x)\ \p_{x_1}\Delta^{-1}\om_k\ \p_{x_2}\Delta^{-1}\om_k\ dx\\[5mm]
\ds=-\int_{{\R}^2}\p_{x_1}\psi(x)\ \ti{\varphi}_k\, \p_{x_2}\ti{\varphi}_k\ dx-\int_{{\R}^2}\psi(x)\ \Delta^{-1}\om_k\ \p^2_{x_1x_2}\Delta^{-1}\om_k\ dx\\[5mm]
\ds=-\int_{{\R}^2}\p_{x_1}\psi(x)\ \ti{\varphi}_k\, \p_{x_2}\ti{\varphi}_k\ dx+\int_{{\R}^2}\lf[\Delta^{-1}(\psi\, \om_k)-\psi(x)\ \Delta^{-1}\om_k\rg]\ \p^2_{x_1x_2}\ti{\varphi}_k\ dx\\[5mm]
\ds\quad+\int_{{\R}^2}\psi(x)\ \om_{k}(x)\ \p^2_{x_1x_2}\Delta^{-2}\om_k(x)\ dx
\end{array}
\ee
We shall now pass to the limit in the three terms in the r.h.s. of (\ref{III.19}). 

\medskip

{\bf The first term of the r.h.s. of (\ref{III.19}).} Since $\ti{\varphi}_k\rightharpoonup\ti{\varphi}_\infty$ weakly in  $W^{1,2}(D^2)$, from Rellich Kondrachoff theorem $\ti{\varphi}_k$ converges
strongly to $\ti{\varphi}_\infty$ in $L^2(D^2)$ therefore
\be
\label{III.20}
\lim_{k\rightarrow +\infty}\int_{{\R}^2}\p_{x_1}\psi(x)\ \ti{\varphi}_k\ \p_{x_2}\ti{\varphi}_k\ dx=\int_{{\R}^2}\p_{x_1}\psi(x)\ \ti{\varphi}_\infty\ \p_{x_2}\ti{\varphi}_\infty\ dx\quad.
\ee

\medskip

{\bf The second term of the r.h.s. of (\ref{III.19}).} Observe first that
\be
\label{III.21}
\Delta\lf[\Delta^{-1}(\psi\, \om_k)-\psi(x)\ \Delta^{-1}\om_k\rg]=-\Delta\psi\ \ti{\varphi}_k-2\nabla\psi\,\nabla \ti{\varphi}_k\quad.
\ee
Since $\ti{\varphi}_k\rightharpoonup\ti{\varphi}_\infty$ weakly in  $W^{1,2}({\R}^2)$, we have that
\be
\label{III.22}
\Delta\lf[\Delta^{-1}(\psi\, \om_k)-\psi(x)\ \Delta^{-1}\om_k\rg]\rightharpoonup \Delta\lf[\Delta^{-1}(\psi\, \om_\infty)-\psi(x)\ \Delta^{-1}\om_\infty\rg]\quad\quad\mbox{ weakly in }\quad L^{2}({\R}^2)
\ee
Hence, using again Rellich-Kondrachoff we deduce that
\be
\label{III.23}
\lf[\Delta^{-1}(\psi\, \om_k)-\psi(x)\ \Delta^{-1}\om_k\rg]\longrightarrow\lf[\Delta^{-1}(\psi\, \om_\infty)-\psi(x)\ \Delta^{-1}\om_\infty\rg]\quad\quad\mbox{ strongly in }\quad W^{1,2}({\R}^2)
\ee
Since
\be
\label{III.24}
\p^2_{x_1x_2}\ti{\varphi}_k\rightharpoonup \p^2_{x_1x_2}\ti{\varphi}_\infty\quad\quad\mbox{ weakly in }\quad H^{-1}({\R}^2)
\ee
Combining (\ref{III.23}) and (\ref{III.24}) gives
\be
\label{III.25}
\lim_{k\rightarrow +\infty}\int_{{\R}^2}\lf[\Delta^{-1}(\psi\, \om_k)-\psi(x)\ \Delta^{-1}\om_k\rg]\ \p^2_{x_1x_2}\ti{\varphi}_k\ dx=\int_{{\R}^2}\lf[\Delta^{-1}(\psi\, \om_\infty)-\psi(x)\ \Delta^{-1}\om_\infty\rg]\ \p^2_{x_1x_2}\ti{\varphi}_\infty\ dx\quad.
\ee

\medskip

{\bf The third term of the r.h.s. of (\ref{III.19}).} This is of course the most delicate term in which the specificity of the bilinearity  $\p_{x_1}\ti{\varphi}_k\,\p_{x_2}\ti{\varphi}_k$ we are considering plays a role.

\medskip

From \cite{Ste} we have that the Kernel associated to the operator  $\p^2_{x_1x_2}\Delta^{-2}$ {\bf is bounded in $L^\infty$}. Indeed one has that the Fourier multiplier associated to
 the operator $\p^2_{x_1x_2}\Delta^{-2}$ is given by
\be
\label{III.26}
\widehat{\p^2_{x_1x_2}\Delta^{-2}}=-\frac{\xi_1\ \xi_2}{|\xi|^4}
\ee
which as to be understood either as in a singular integral sense or in distributional sense as being the following tempered distribution in ${\mathcal S}'({\R}^2)$
\[
\begin{array}{l}
\ds\forall \ \phi(\xi)\in {\mathcal S}({\R}^2)\quad\quad\lf<pv\lf(-\frac{\xi_1\ \xi_2}{|\xi|^4}\rg);\phi(\xi)\rg>=-\lim_{\ep\rightarrow 0}\int_{{\R}^2\setminus B_\ep(0)}\frac{\xi_1\ \xi_2}{|\xi|^4}\ \phi(\xi)\ d\xi\\[5mm]
\ds\quad\quad\quad=
\int_{{\R}^2}\frac{\xi_1\ \xi_2}{|\xi|^4}\ (\phi(0)-\phi(\xi))\ d\xi
\end{array}
\]
Since the homogeneous polynomial $\xi_1\, \xi_2$ is harmonic we can apply theorem 5 in 3.3 of \cite{Ste} and deduce the existence of a universal constant $c_0$ such that
the inverse of the Fourier transform of ${\xi_1\ \xi_2}/{|\xi|^4}$ is given by
\[
-\widehat{\frac{\xi_1\ \xi_2}{|\xi|^4}}^{-1}=c_0\ \frac{x_1\ x_2}{|x|^2}\quad.
\]
Hence
\be
\label{III.27}
\int_{{\R}^2}\psi(x)\ \om_{k}(x)\ \p^2_{x_1x_2}\Delta^{-2}\om_k(x)\ dx=c_0\int_{{\R}^2}\psi(x)\ \om_{k}(x)\ \om_k(y)\ \frac{(x_1-y_1)\ (x_2-y_2)}{|x-y|^2}\ dx\ dy
\ee
If the kernel ${(x_1-y_1)\ (x_2-y_2)}/{|x-y|^2}$ would have been continuous up to the diagonal $x=y$ (or even VMO on ${\R}^4$ ) we could have easily pass to the limit in this integral,
since $\om_{k}(x)\ \om_k(y)$ is uniformly bounded in the local Hardy space ${\mathcal H}^1_{loc}(R^4)$, it converges weakly in particular in Radon measure 
to $\om_\infty(x)\ \om_\infty(y)$. We shall however make use of the fact that ${(x_1-y_1)\ (x_2-y_2)}/{|x-y|^2}$ is bounded in $L^\infty$ in order to pass to the limit in (\ref{III.27}) modulo
possible concentration points. 

\medskip

Let $\chi$ be a cut-off function in $C^\infty_0({\R}^+,{\R}^+)$ such that $\chi$ is equal to $1$ on $[0,1]$ and equal to zero on $[2,+\infty)$ and $0\le\chi\le 1$. For $\ep>0$
we denote $\chi_\ep(t):=\chi(t/\ep)$.

\medskip

We write
\be
\label{III.28}
\begin{array}{l}
\ds\int_{{\R}^4}\psi(x)\ \om_{k}(x)\ \om_k(y)\ \frac{(x_1-y_1)\ (x_2-y_2)}{|x-y|^2}\ dx\ dy\\[5mm]
\ds\quad\quad\quad=\int_{{\R}^4}\psi(x)\ \om_{k}(x)\ \om_k(y)\ \chi_\ep(|x-y|)\ \frac{(x_1-y_1)\ (x_2-y_2)}{|x-y|^2}\ dx\ dy\\[5mm]
\ds\quad\quad\quad+\int_{{\R}^4}\psi(x)\ \om_{k}(x)\ \om_k(y)\ [1-\chi_\ep(|x-y|)]\ \frac{(x_1-y_1)\ (x_2-y_2)}{|x-y|^2}\ dx\ dy\quad.
\end{array}
\ee
Since $[1-\chi_\ep(|x-y|)]\ {(x_1-y_1)/ (x_2-y_2)}/{|x-y|^2}$ is continuous on ${\R}^4$ we have
\be
\label{III.29}
\begin{array}{l}
\ds\lim_{k\rightarrow+\infty}\int_{{\R}^4}\psi(x)\ \om_{k}(x)\ \om_k(y)\ [1-\chi_\ep(|x-y|)]\ \frac{(x_1-y_1)\ (x_2-y_2)}{|x-y|^2}\ dx\ dy\\[5mm]
\ds\quad\quad\quad=\int_{{\R}^4}\psi(x)\ \om_{\infty}(x)\ \om_{\infty}(y)\ [1-\chi_\ep(|x-y|)]\ \frac{(x_1-y_1)\ (x_2-y_2)}{|x-y|^2}\ dx\ dy
\end{array}
\ee
And then
\be
\label{III.30}
\begin{array}{l}
\ds\lim_{\ep\rightarrow 0}\lim_{k\rightarrow+\infty}\int_{{\R}^4}\psi(x)\ \om_{k}(x)\ \om_k(y)\ [1-\chi_\ep(|x-y|)]\ \frac{(x_1-y_1)\ (x_2-y_2)}{|x-y|^2}\ dx\ dy\\[5mm]
\ds\quad\quad\quad=\int_{{\R}^4}\psi(x)\ \om_{\infty}(x)\ \om_{\infty}(y)\ \frac{(x_1-y_1)\ (x_2-y_2)}{|x-y|^2}\ dx\ dy
\end{array}
\ee
Combining (\ref{III.20}), (\ref{III.25}) and (\ref{III.30}) we obtain that
\be
\label{III.31}
\lf|\lf<\gamma-\p_{x_1}\ti{\varphi}_\infty\,\p_{x_2}\ti{\varphi}_\infty\ ;\ \psi\rg>\rg|\le \liminf_{\ep\rightarrow 0}\liminf_{k\rightarrow+\infty}\int_{{\R}^4}|\psi(x)|\ |\om_{k}(x)|\ |\om_k(y)|\ \chi_\ep(|x-y|)\ dx\ dy
\ee
Modulo extraction of a subsequence we can assume that the sequence of measures $|\om_k(x)|\ dx$ converges weakly to a non negative Radon measure $\zeta(x)$ and we have
\be
\label{III.32}
\lf|\lf<\gamma-\p_{x_1}\ti{\varphi}_\infty\,\p_{x_2}\ti{\varphi}_\infty\ ;\ \psi\rg>\rg|\le \liminf_{\ep\rightarrow 0}<|\psi(x)|\ \chi_\ep(|x-y|)\, ;\ \zeta(x)\ \zeta(y)>
\ee
Denote by $A_\zeta:=\sum_{j\in {\N}}\zeta_j\ \delta_{q_j}$  the atomic part of $\zeta$ : 
\[
<\zeta(y)-A_\zeta(y)\ ;\ \chi_\ep(|x-y|)>\rightarrow 0\quad\quad\zeta\mbox{ a. e. }x\quad.
\]
Thus
\[
\lim_{\ep\rightarrow 0}<|\psi(x)|\ \chi_\ep(|x-y|)\, ;\ \zeta(x)\ \zeta(y)>=\lim_{\ep\rightarrow 0}\sum_{j\in {\N}}\ \zeta_j\ <|\psi(x)|\ \chi_\ep(|x-q_j|)\, ;\ \zeta(x)>=\sum_{j\in {\N}}\ \zeta_j^2\ |\psi(q_j)|
\]
Hence (\ref{III.32}) implies
\be
\label{III.33}
\lf|\lf<\gamma-\p_{x_1}\ti{\varphi}_\infty\,\p_{x_2}\ti{\varphi}_\infty\ ;\ \psi\rg>\rg|\le \sum_{j\in {\N}}\ \zeta_j^2\ |\psi(q_j)|
\ee
which shows that $\gamma-\p_{x_1}\ti{\varphi}_\infty\,\p_{x_2}\ti{\varphi}_\infty$ is atomic. This implies the lemma for the bilinearity $\p_{x_1}\varphi_\infty\,\p_{x_2}\varphi_\infty$.
The same applies to the bilinearity $|\p_{x_1}\varphi_\infty|^2-|\p_{x_2}\varphi_\infty|^2$ since in the estimation of the third term in the r.h.s of the identity corresponding to (\ref{III.19})
one uses that $\xi_1^2-\xi_2^2$ is harmonic and thus the kernel associated to $(\p_{x_1^2}^2-\p_{x_2^2}^2)\Delta^{-2}$ is also bounded in $L^\infty$ due to theorem 5 in section 3.3 of \cite{Ste}.\hfill $\Box$

\section{Proof of the main theorem~\ref{th-I.3}.}
\reset
Let $\vec{\Phi}_k$ be a sequence of global weak isothermic immersions of an abstract closed surface $\Sigma$ into ${\R}^m$ such that the conformal class  to which the induced
metric $g_k:=\vec{\Phi}_k^\ast g_{{\R}^m}$ does not degenerate. This means that there exists a sequence of constant scalar curvature metric $h_k$ of volume 1,
precompact for any $C^l$ norm of $\Sigma$ (equipped with some fixed arbitrary reference metric $g_0$) and a diffeomorphism $f_k$ of $\Sigma$ such that
\be
\label{III.1}
\vec{\Phi}_k\circ f_k\quad :\quad (\Sigma,h_k)\longrightarrow {\R}^m\quad\quad\mbox{ is conformal }\quad.
\ee
Modulo extraction of a subsequence we can assume that
\be
\label{III.2}
h_k\rightarrow h_\infty\quad\quad\mbox{ in }C^l(\Sigma)\quad\forall l\in {\N}\quad,
\ee
where $h_\infty$ is a constant scalar curvature of volume 1 on $\Sigma$.

We assume moreover that
\[
\limsup_{k\rightarrow +\infty}W(\vec{\Phi}_k)=\limsup_{k\rightarrow+\infty}\int_\Sigma|\vec{H}_{\vec{\Phi}_k}|^2\ dvol_{g_k}<+\infty
\]
Following the normalization lemma A.4 and lemma III.1 of \cite{Ri3}, we deduce the existence of a sequence of M\"obius transformation $\Xi_k$ of ${\R}^m\cup\{\infty\}$
(i.e. $\Xi_k$ are conformal diffeomorphism of ${\R}^m\cup\{\infty\}$) such that $\vec{\xi}_k:=\Xi_k\circ\vec{\Phi}_k\circ f_k$ satisfies the following conditions (up to subsequence)
\begin{itemize}
\item[i)] 
\[
\exists \, R>0\quad,\quad\forall\, k\quad\quad \vec{\xi}_k(\Sigma)\subset B_R(0)\quad.
\]
\item[ii)]
\[
\exists \ a_1\cdots a_N\quad\quad\mbox{ s.t. }\quad\quad\vec{\xi}_k\rightharpoonup \vec{\xi}_\infty\quad\quad\mbox{ in }W^{2,2}_{loc}(\Sigma\setminus\{a_1\cdots a_N\})
\]
\item[iii)]
\[
\forall\, K\ \mbox{ compact of }\Sigma\setminus\{a_1\cdots a_N\}\quad\quad\limsup_{k\rightarrow +\infty}\|\log|d\vec{\xi}_k|_{h_k}\|_{L^\infty(K)}<+\infty
\]
\end{itemize}

These 3 conditions ensure that the weak limiting map $\vec{\xi}_\infty$ is a weak possibly branched conformal  immersion in the space ${\mathcal F}_\Sigma$.

\medskip

Assuming now that $\vec{\Phi}_k$ are {\it weak global isothermic immersions} in ${\mathcal E}_\Sigma$ then, due to the conformal invariance proved in proposition~\ref{pr-I.4},
$\vec{\xi}_k$ are also {\it weak global isothermic immersions}. Thus there exists a sequence of non zero holomorphic quadratic differentials $q_k$ for the sequence of riemann surfaces
$(\Sigma,h_k)$ satisfying
\be
\label{III.3}
\Im(q_k,\vec{h}_{0,k})_{WP}=0
\ee
where the Weil-Peterson norm is taken with respect to $h_k$. Because of the linearity of equation (\ref{III.3}) with respect to $q_k$ we can normalize $q_k$ in such a way that
\be
\label{III.4}
\forall k\in{\N}\quad\quad\quad\int_{\Sigma}(q_k,q_k)_{WP}\ dvol_{{h}_k}=1
 \ee
 The space $P_k$ of holomorphic quadratic forms of $(\Sigma,h_k)$ is a finite dimensional space of fixed dimension (depending on $\Sigma$ only) of the space $\Gamma(T^\ast\Sigma\otimes T^\ast\Sigma)$  of smooth sections of $T^\ast\Sigma\otimes T^\ast\Sigma$. Since $h_k$ converges to $h_\infty$ we can extract a subsequence such that $P_k$ converges to $P_\infty$
 and we can extract a subsequence such that $q_k$ converges in any $C^l$ norm towards $q_\infty$ for any $l\in {\N}$.
 
 \medskip
 
 The holomorphic quadratic form $q_\infty$ satisfy also (\ref{III.4}), moreover, due to the weak convergence of $\vec{h}_{0,k}$ towards $\vec{h}_{0,\infty}$ in $L^2_{loc}(\Sigma\setminus\{a_1\cdots a_N\})$, 
\be
\label{III.5}
\Im(q_\infty,\vec{h}_{0,\infty})_{WP}=0\quad\quad\quad\mbox{ in }\Sigma\setminus\{a_1\cdots a_N\}\quad.
\ee
This implies that $\vec{\xi}_\infty$ is a weak, possibly branched, conformal isothermic immersion of $(\Sigma,h_\infty)$ into ${\R}^m$.

\medskip

In an arbitrary strongly converging conformal chart $\phi_k \ :\ D^2\setminus (\Sigma,h_k)$ the equation satisfied by $\vec{\xi}_k\circ\phi_k$ reads (omitting to write explicitly the composition with 
$\phi_k$)
\[
\Im(f_k(z)\ \ov{\vec{H}_{0,k}})=2\,\Im\lf(f_k(z)\ \p_{\ov{z}}\lf[e^{-2\,\la_k}\ \partial_{\ov{z}}\vec{\xi}_k\rg]\rg)
\]
where $f_k$ is the expression of $q_k$ in this chart $f_k(z)\ dz\otimes dz=q_k$.

\medskip

Denote by $b_1\cdots b_Q$ the isolated zeros of $q_\infty$ in $\Sigma$. Let $U$ be a disc included in $\Sigma\setminus\{a_1\cdots a_N,b_1\cdots b_Q\}$. Considering a converging sequence of conformal charts $\phi_k$ realizing a diffeomorphism from $D^2$ into $U$, since $f_\infty$ the expression of $\vec{h}_{0,\infty}$ in this chart does not vanish on $D^2$ and since $f_k$ converge strongly on $D^2$ towards $f_\infty$, we can introduce the new {\bf converging} chart $w:=\sqrt{f_k\circ\phi_k^{-1}}$. In this new chart the isothermic equation reads
\be
\label{III.6}
\Im(\vec{H}_{0,k})=\frac{\p}{\p x_1}\lf[e^{-2\la_k}\frac{\p\vec{\xi}_k}{\p x_2}\rg]+\frac{\p}{\p x_2}\lf[e^{-2\la_k}\frac{\p\vec{\xi}_k}{\p x_1}\rg]=0\quad.
\ee
where $w=x_1+ix_2$ and $e^{\la_k}=|\p_{x_1}\vec{\xi}_k|=|\p_{x_2}\vec{\xi}_k|$. since the chart is strongly converging the expression of $\vec{\xi}_k$ in this chart satisfy
\be
\label{III.7}
\vec{\xi}_k(w)\rightharpoonup\vec{\xi}_\infty(w)\quad\mbox{ in }\quad W^{2,2}(D^2)\quad\mbox{  and }\quad\quad\limsup_{k\rightarrow +\infty}\|\la_k(w)\|_{L^\infty(D^2)}<+\infty\quad.
\ee
We also choose $U$ small enough and the subsequence in such a way that
\be
\label{III.8}
\forall\, k\in {\N}\quad\int_{D^2}|\nabla \vec{n}_{\vec{\xi}_k}|^2\ dx_1\, dx_2<\frac{8\pi}{3}\quad.
\ee
We can then use a result by F. H\'elein (see \cite{He} chapter 5) that gives the existence of $(\vec{e}_{1,k},\vec{e}_{2,k})\in (W^{1,2}(D^2,S^{m-1}))^2$ such that
\be
\label{III.9}
\vec{e}_{1,k}\wedge\vec{e}_{2,k}=\star\,\vec{n}_{\vec{\xi}_k}\quad\quad\int_{D^2}\sum_{i=1}^2|\nabla\vec{e}_{i,k}|^2<C\,\int_{D^2}|\nabla\vec{n}_{\vec{\xi}_k}|^2
\ee
where $C$ is independent of $k$. We can use this moving frame to express the laplacian of $\la_k$ (see \cite{Ri1}) and we have precisely
\be
\label{III.10}
-\Delta\la_k=\p_{x_1}\vec{e}_{1,k}\cdot\p_{x_2}\vec{e}_{2,k}-\p_{x_2}\vec{e}_{1,k}\cdot\p_{x_1}\vec{e}_{2,k}\quad\quad\mbox{ in }D^2\quad.
\ee
Let $s_k$ be the solution of
\be
\label{III.11a}
\lf\{
\begin{array}{l}
\ds-\Delta s_k=\p_{x_1}\vec{e}_{1,k}\cdot\p_{x_2}\vec{e}_{2,k}-\p_{x_2}\vec{e}_{1,k}\cdot\p_{x_1}\vec{e}_{2,k}\quad\quad\mbox{ in }D^2\\[5mm]
\ds \quad s_k=0\quad\quad\quad\mbox{ on }\quad\partial D^2
\end{array}
\rg.
\ee
From Wente theorem (see \cite{We} and \cite{He}) we have
\be
\label{III.12a}
\|s_k\|_{L^\infty(D^2)}\le C\ \|\nabla\vec{e}_{1,k}\|_{L^2(D^2)}\ \|\nabla\vec{e}_{2,k}\|_{L^2(D^2)}\le C\,\int_{D^2}|\nabla\vec{n}_{\vec{\xi}_k}|^2\quad.
\ee
Using (\ref{III.8}) we deduce that $s_k$ is uniformly bounded in $L^\infty(D^2)$. Combining this fact with (\ref{III.7}) we obtain that the harmonic function
$v_k:=\la_k-s_k$ is uniformly  bounded in $L^\infty(D^2)$. Thus we have that
\be
\label{III.13a}
v_k\rightarrow v_\infty\quad\quad\quad\mbox{ in }\quad C^l_{loc}(D^2)\quad\forall l\in{\N}\quad.
\ee
Lemma~\ref{lm-III.1} implies that there exists a subsequence and two atomic measures $\mu$ and $\nu$ such that
there exists $p_i\in D^2$ for $i\in {\N}$, and $q_j\in D^2$ for $j\in {\N}$ satisfying
\be
\label{III.15a}
\mu=\sum_{i\in {\N}}c_i\, \delta_{p_i}\quad\quad\mbox{ and }\quad\quad\nu=\sum_{j\in {\N}}d_j\, \delta_{q_j}\quad\quad\mbox{in }{\mathcal D}'(D^2)\quad,
\ee
where
\be
\label{III.16a}
\sum_{i\in {\N}}|c_i|=|\mu|(D^2)<+\infty\quad\quad\mbox{ and }\quad\quad\sum_{j\in {\N}}|d_j|=|\nu|(D^2)<+\infty\quad.
\ee
and
\be
\label{III.17z}
\lf\{
\begin{array}{l}
|\p_{x_1}s_{k'}|^2-|\p_{x_2}s_{k'}|^2\rightharpoonup|\p_{x_1}s_{\infty}|^2-|\p_{x_2}s_{\infty}|^2+\mu\quad\quad\mbox{ in }\quad{\mathcal D}'(D^2)\\[5mm]
\p_{x_1}s_{k'}\ \p_{x_2}s_{k'}\rightharpoonup\p_{x_1}s_{\infty}\ \p_{x_2}s_{\infty}+\nu\quad\quad\mbox{ in }\quad{\mathcal D}'(D^2)\quad .
\end{array}
\rg.
\ee
Using (\ref{III.13a}) we deduce 
\be
\label{III.18z}
\lf\{
\begin{array}{l}
|\p_{x_1}\la_{k'}|^2-|\p_{x_2}\la_{k'}|^2\rightharpoonup|\p_{x_1}\la_{\infty}|^2-|\p_{x_2}\la_{\infty}|^2+\mu\quad\quad\mbox{ in }\quad{\mathcal D}'(D^2)\\[5mm]
\p_{x_1}\la_{k'}\ \p_{x_2}\la_{k'}\rightharpoonup\p_{x_1}\la_{\infty}\ \p_{x_2}\la_{\infty}+\nu\quad\quad\mbox{ in }\quad{\mathcal D}'(D^2)\quad .
\end{array}
\rg.
\ee
Assuming the $\vec{\xi}_{k'}$ and  $\vec{\xi}_\infty$ are smooth, since these immersions are smooth, we can apply proposition~\ref{pr-II.1} and deduce that in one hand
\be
\label{III.19z}
\lf\{
\begin{array}{l}
\ds\frac{\p}{\p x_1}\lf[ \lf(\frac{\p\vec{n}_{\vec{\xi}_{k'}}}{\p x_2}\res\vec{e}_{2,k'}\rg)^2+\lf|\frac{\p\la_{k'}}{\p{x_1}}\rg|^2-\lf|\frac{\p\la_{k'}}{\p{x_2}}\rg|^2\rg]+\frac{\p}{\p {x_2}}\lf[2\ \frac{\p\la_{k'}}{\p{x_1}}\ \frac{\p\la_{k'}}{\p{x_2}}\rg]=0\\[5mm]
\ds\frac{\p}{\p x_2}\lf[ \lf(\frac{\p\vec{n}_{\vec{\xi}_{k'}}}{\p x_1}\res\vec{e}_{1,k'}\rg)^2+\lf|\frac{\p\la_{k'}}{\p{x_2}}\rg|^2-\lf|\frac{\p\la_{k'}}{\p{x_1}}\rg|^2\rg]+\frac{\p}{\p {x_1}}\lf[2\ \frac{\p\la_{k'}}{\p{x_1}}\ \frac{\p\la_{k'}}{\p{x_2}}\rg]=0
\end{array}
\rg.
\ee
and in the other hand
\be
\label{III.20z}
\lf\{
\begin{array}{l}
\ds\frac{\p}{\p x_1}\lf[ \lf(\frac{\p\vec{n}_{\vec{\xi}_{\infty}}}{\p x_2}\res\vec{e}_{2,\infty}\rg)^2+\lf|\frac{\p\la_{\infty}}{\p{x_1}}\rg|^2-\lf|\frac{\p\la_{\infty}}{\p{x_2}}\rg|^2\rg]+\frac{\p}{\p {x_2}}\lf[2\ \frac{\p\la_{\infty}}{\p{x_1}}\ \frac{\p\la_{\infty}}{\p{x_2}}\rg]=0\\[5mm]
\ds\frac{\p}{\p x_2}\lf[ \lf(\frac{\p\vec{n}_{\vec{\xi}_{\infty}}}{\p x_1}\res\vec{e}_{1,\infty}\rg)^2+\lf|\frac{\p\la_{\infty}}{\p{x_2}}\rg|^2-\lf|\frac{\p\la_{\infty}}{\p{x_1}}\rg|^2\rg]+\frac{\p}{\p {x_1}}\lf[2\ \frac{\p\la_{\infty}}{\p{x_1}}\ \frac{\p\la_{\infty}}{\p{x_2}}\rg]=0
\end{array}
\rg.
\ee
Applying Poincar\'e Lemma, we deduce the existence of $A_{k'}$ and  $B_{k'}$ in $W^{1,1}$ such that
\[
\lf\{
\begin{array}{l}
\ds\p_{x_2}A_{k'}= \lf(\frac{\p\vec{n}_{\vec{\xi}_{k'}}}{\p x_2}\res\vec{e}_{2,k'}\rg)^2+\lf|\frac{\p\la_{k'}}{\p{x_1}}\rg|^2-\lf|\frac{\p\la_{k'}}{\p{x_2}}\rg|^2\\[5mm]
\ds\p_{x_1} A_{k'}=-2\ \frac{\p\la_{k'}}{\p{x_1}}\ \frac{\p\la_{k'}}{\p{x_2}}
\end{array}
\rg.
\]
such that
\[
\lf\{
\begin{array}{l}
\ds\p_{x_1}B_{k'}=\lf(\frac{\p\vec{n}_{\vec{\xi}_{k'}}}{\p x_1}\res\vec{e}_{1,k'}\rg)^2+\lf|\frac{\p\la_{k'}}{\p{x_2}}\rg|^2-\lf|\frac{\p\la_{k'}}{\p{x_1}}\rg|^2\\[5mm]
\ds\p_{x_2} B_{k'}=-2\ \frac{\p\la_{k'}}{\p{x_1}}\ \frac{\p\la_{k'}}{\p{x_2}}
\end{array}
\rg.
\]
Moreover for the same reason there exist $A_{\infty}$ and  $B_{\infty}$ in $W^{1,1}$ such that
\[
\lf\{
\begin{array}{l}
\ds\p_{x_2}A_{\infty}= \lf(\frac{\p\vec{n}_{\vec{\xi}_{\infty}}}{\p x_2}\res\vec{e}_{2,\infty}\rg)^2+\lf|\frac{\p\la_{\infty}}{\p{x_1}}\rg|^2-\lf|\frac{\p\la_{\infty}}{\p{x_2}}\rg|^2\\[5mm]
\ds\p_{x_1} A_{\infty}=-2\ \frac{\p\la_{\infty}}{\p{x_1}}\ \frac{\p\la_{\infty}}{\p{x_2}}
\end{array}
\rg.
\]
such that
\[
\lf\{
\begin{array}{l}
\ds\p_{x_1}B_{\infty}=\lf(\frac{\p\vec{n}_{\vec{\xi}_{\infty}}}{\p x_1}\res\vec{e}_{1,\infty}\rg)^2+\lf|\frac{\p\la_{\infty}}{\p{x_2}}\rg|^2-\lf|\frac{\p\la_{\infty}}{\p{x_1}}\rg|^2\\[5mm]
\ds\p_{x_2} B_{\infty}=-2\ \frac{\p\la_{\infty}}{\p{x_1}}\ \frac{\p\la_{\infty}}{\p{x_2}}
\end{array}
\rg.
\]
We observe that we have
\[
\p_{x_1}A_{k'}=\p_{x_2}B_{k'}\quad\quad\mbox{ and }\quad\quad\p_{x_1}A_{\infty}=\p_{x_2}B_{\infty}
\]
Applying again Poincar\'e Lemma, we have the existence of $\al_{k'}$ and $\al_\infty$ in $W^{2,1}$ such that
\[
\nabla\al_{k'}=(B_{k'},A_{k'})\quad\quad\quad\mbox{ and }\quad\quad\quad\nabla\al_{\infty}=(B_{\infty},A_{\infty})
\]
Thus we have
\be
\label{III.20a}
\lf\{
\begin{array}{l}
\ds\frac{\p^2\al_{k'}}{\p x_2\p x_2}= \lf(\frac{\p\vec{n}_{\vec{\xi}_{k'}}}{\p x_2}\res\vec{e}_{2,k'}\rg)^2+\lf|\frac{\p\la_{k'}}{\p{x_1}}\rg|^2-\lf|\frac{\p\la_{k'}}{\p{x_2}}\rg|^2\\[5mm]
\ds\frac{\p^2\al_{k'}}{\p x_1\p x_1}=\lf(\frac{\p\vec{n}_{\vec{\xi}_{k'}}}{\p x_1}\res\vec{e}_{1,k'}\rg)^2+\lf|\frac{\p\la_{k'}}{\p{x_2}}\rg|^2-\lf|\frac{\p\la_{k'}}{\p{x_1}}\rg|^2\\[5mm]
\ds\frac{\p^2\al_{k'}}{\p x_1\p x_2}=-2\ \frac{\p\la_{k'}}{\p{x_1}}\ \frac{\p\la_{k'}}{\p{x_2}}
\end{array}
\rg.
\ee
and
\be
\label{III.20b}
\lf\{
\begin{array}{l}
\ds\frac{\p^2\al_{\infty}}{\p x_2\p x_2}= \lf(\frac{\p\vec{n}_{\vec{\xi}_{\infty}}}{\p x_2}\res\vec{e}_{2,\infty}\rg)^2+\lf|\frac{\p\la_{\infty}}{\p{x_1}}\rg|^2-\lf|\frac{\p\la_{\infty}}{\p{x_2}}\rg|^2\\[5mm]
\ds\frac{\p^2\al_{\infty}}{\p x_1\p x_1}=\lf(\frac{\p\vec{n}_{\vec{\xi}_{\infty}}}{\p x_1}\res\vec{e}_{1,\infty}\rg)^2+\lf|\frac{\p\la_{\infty}}{\p{x_2}}\rg|^2-\lf|\frac{\p\la_{\infty}}{\p{x_1}}\rg|^2\\[5mm]
\ds\frac{\p^2\al_{\infty}}{\p x_1\p x_2}=-2\ \frac{\p\la_{\infty}}{\p{x_1}}\ \frac{\p\la_{\infty}}{\p{x_2}}
\end{array}
\rg.
\ee
Since $\nabla A_{k'}$ and $\nabla B_{k'}$ are uniformly bounded in $L^1$ we can normalize $A_{k'}$ and $B_{k'}$ in such a way that
$A_{k'}$ and $B_{k'}$ are uniformly bounded in $L^2(D^2)$. In a similar way, since now $\nabla\al_{k'}$ is uniformly bounded 
in $L^2(D^2)$ we can normalize $\al_{k'}$ in such a way that $\al_{k'}$ is uniformly bounded in $W^{1,2}(D^2)$. Thus
\be
\label{III.21z}
\limsup_{k'\rightarrow +\infty}\lf\|\lf(\frac{\p\vec{n}_{\vec{\xi}_{k'}}}{\p x_2}\res\vec{e}_{2,k'}\rg)^2+\lf|\frac{\p\la_{k'}}{\p{x_1}}\rg|^2-\lf|\frac{\p\la_{k'}}{\p{x_2}}\rg|^2\rg\|_{H^{-1}(D^2)}<+\infty\quad,
\ee
moreover
\be
\label{III.22z}
\limsup_{k'\rightarrow +\infty}\lf\|\lf(\frac{\p\vec{n}_{\vec{\xi}_{k'}}}{\p x_1}\res\vec{e}_{1,k'}\rg)^2+\lf|\frac{\p\la_{k'}}{\p{x_2}}\rg|^2-\lf|\frac{\p\la_{k'}}{\p{x_1}}\rg|^2\rg\|_{H^{-1}(D^2)}<+\infty\quad,
\ee
and finally
\be
\label{III.23z}
\limsup_{k'\rightarrow +\infty}\lf\|\frac{\p\la_{k'}}{\p{x_1}}\ \frac{\p\la_{k'}}{\p{x_2}}\rg\|_{H^{-1}(D^2)}<+\infty\quad.
\ee
Taking this last quantity, we can always extract a subsequence, that we will still denote $k'$, such that
\be
\label{III.24z}
\frac{\p\la_{k'}}{\p{x_1}}\ \frac{\p\la_{k'}}{\p{x_2}}\rightharpoonup f\quad\quad\mbox{ weakly in }H^{-1}(D^2)
\ee
Comparing this convergence with the second line of (\ref{III.18}) gives
\be
\label{III.25z}
f=\p_{x_1}\la_{\infty}\ \p_{x_2}\la_{\infty}+\sum_{j\in {\N}}d_j\, \delta_{q_j}\in H^{-1}(D^2)\quad.
\ee
But, argueing as for $\al_{k'}$, we have that $\al_{\infty}\in W^{1,1}(D^2)$ and hence, using the last line of (\ref{III.20b}), we have that
\[
\nu=\sum_{j\in {\N}}d_j\, \delta_{q_j}\in H^{-1}(D^2)
\] 
This implies that this atomic measure is zero,
\be
\label{III.26z}
\nu\equiv 0
\ee 
which is the unique atomic measure included in $H^{-1}$.

\medskip

Similarly, from (\ref{III.21z}) we can extract a subsequence, still denoted $k'$, such that
\be
\label{III.27z}
\lf(\frac{\p\vec{n}_{\vec{\xi}_{k'}}}{\p x_2}\res\vec{e}_{2,k'}\rg)^2+\lf|\frac{\p\la_{k'}}{\p{x_1}}\rg|^2-\lf|\frac{\p\la_{k'}}{\p{x_2}}\rg|^2\rightharpoonup g_1\quad\quad\mbox{ weakly in }H^{-1}(D^2)
\ee
and
\be
\label{III.28z}
\lf(\frac{\p\vec{n}_{\vec{\xi}_{k'}}}{\p x_1}\res\vec{e}_{1,k'}\rg)^2+\lf|\frac{\p\la_{k'}}{\p{x_2}}\rg|^2-\lf|\frac{\p\la_{k'}}{\p{x_1}}\rg|^2\rightharpoonup g_1\quad\quad\mbox{ weakly in }H^{-1}(D^2)
\ee
Comparing these convergences with the first line of (\ref{III.18z}) gives in one hand
\be
\label{III.29z}
h_-:=\lf|\frac{\p\la_{\infty}}{\p{x_1}}\rg|^2-\lf|\frac{\p\la_{\infty}}{\p{x_2}}\rg|^2+\sum_{i\in {\N}}c_i\, \delta_{p_i}-g_1\le 0
\ee
and in the other hand
\be
\label{III.30z}
h_+:=\lf|\frac{\p\la_{\infty}}{\p{x_1}}\rg|^2-\lf|\frac{\p\la_{\infty}}{\p{x_2}}\rg|^2+\sum_{i\in {\N}}c_i\, \delta_{p_i}+g_2\ge 0
\ee
Using the two first lines of (\ref{III.20b}), we have that
\[
\lf[\lf|\frac{\p\la_{\infty}}{\p{x_1}}\rg|^2-\lf|\frac{\p\la_{\infty}}{\p{x_2}}\rg|^2-g_1\rg]\in  H^{-1}(D^2)\quad\quad\mbox{ and }\quad\quad\lf[\lf|\frac{\p\la_{\infty}}{\p{x_1}}\rg|^2-\lf|\frac{\p\la_{\infty}}{\p{x_2}}\rg|^2+g_2
\rg]\in  H^{-1}(D^2)
\]
Let $\chi$ be a cut off function in $C^\infty({\R}^+,{\R}^+)$ identically equal to 1 on $[0,1]$, equal to $0$ on $(2,+\infty)$ and $0\le \chi\le 1$ on ${\R}^+$.  For any $\ep>0$ we denote
$\chi_\ep(t):=\chi(t/\ep)$. For any $i_0\in {\N}$, the map $\chi_\ep(|x-p_{i_0}|)$ weakly converge to zero in $W^{1,2}_0(D^2)$ thus
\be
\label{III.31z}
0\ge\lim_{\ep\rightarrow 0}\lf<h_-,\chi_\ep(|x-p_{i_0}|)\rg>=\lim_{\ep\rightarrow 0}\lf<\lf|\frac{\p\la_{\infty}}{\p{x_1}}\rg|^2-\lf|\frac{\p\la_{\infty}}{\p{x_2}}\rg|^2+\sum_{i\in {\N}}c_i\, \delta_{p_i}-g_1,\chi_\ep(|x-p_{i_0}|)\rg>=c_{i_0}
\ee
In a similar way we have
\be
\label{III.32z}
0\le\lim_{\ep\rightarrow 0}<h_+,\chi_\ep(|x-p_{i_0}|)>=\lim_{\ep\rightarrow 0}\lf<\lf|\frac{\p\la_{\infty}}{\p{x_1}}\rg|^2-\lf|\frac{\p\la_{\infty}}{\p{x_2}}\rg|^2+\sum_{i\in {\N}}c_i\, \delta_{p_i}+g_2,\chi_\ep(|x-p_{i_0}|)\rg>
= c_{i_0}
\ee
Comparing (\ref{III.31z}) and (\ref{III.32z}) gives for any $i_0$ $c_{i_0}=0$ and then we have proved that
\be
\label{III.33z}
\mu\equiv 0
\ee 
Combining (\ref{III.18z}), (\ref{III.26z}) and (\ref{III.32z}) implies then
\be
\label{III.34z}
\lf\{
\begin{array}{l}
|\p_{x_1}\la_{k'}|^2-|\p_{x_2}\la_{k'}|^2\rightharpoonup|\p_{x_1}\la_{\infty}|^2-|\p_{x_2}\la_{\infty}|^2\quad\quad\mbox{ in }\quad{\mathcal D}'(D^2)\\[5mm]
\p_{x_1}\la_{k'}\ \p_{x_2}\la_{k'}\rightharpoonup\p_{x_1}\la_{\infty}\ \p_{x_2}\la_{\infty}\quad\quad\mbox{ in }\quad{\mathcal D}'(D^2)\quad .
\end{array}
\rg.
\ee
Translating this information in terms of $\al_{k'}$ and $\al_\infty$ gives
\be
\label{III.35z}
\lf\{
\begin{array}{l}
\ds-\frac{\p^2\al_{k'}}{\p x_2\p x_2}+ \lf(\frac{\p\vec{n}_{\vec{\xi}_{k'}}}{\p x_2}\res\vec{e}_{2,k'}\rg)^2\rightharpoonup -\frac{\p^2\al_{\infty}}{\p x_2\p x_2}+ \lf(\frac{\p\vec{n}_{\vec{\xi}_{\infty}}}{\p x_2}\res\vec{e}_{2,\infty}\rg)^2\quad\quad\mbox{ weakly in }H^{-1}(D^2)\\[5mm]
\ds-\frac{\p^2\al_{k'}}{\p x_1\p x_1}+\lf(\frac{\p\vec{n}_{\vec{\xi}_{k'}}}{\p x_1}\res\vec{e}_{1,k'}\rg)^2\rightharpoonup -\frac{\p^2\al_{\infty}}{\p x_1\p x_1}+\lf(\frac{\p\vec{n}_{\vec{\xi}_{\infty}}}{\p x_1}\res\vec{e}_{1,\infty}\rg)^2\quad\quad\mbox{ weakly in }H^{-1}(D^2)\\[5mm]
\ds\frac{\p^2\al_{k'}}{\p x_1\p x_2}\rightharpoonup\frac{\p^2\al_{\infty}}{\p x_1\p x_2}\quad\quad\mbox{ weakly in }H^{-1}(D^2)
\end{array}
\rg.
\ee
Denote by $\hat{\al}_\infty$ the weak limit (modulo extraction of a subsequence) of $\al_{k'}$ in $W^{1,2}$ and let $\beta_\infty:=\hat{\al}_\infty-\al_\infty$. We have
\be
\label{III.36z}
\lf\{
\begin{array}{l}
\ds \lf(\frac{\p\vec{n}_{\vec{\xi}_{k'}}}{\p x_1}\res\vec{e}_{1,k'}\rg)^2+\lf(\frac{\p\vec{n}_{\vec{\xi}_{k'}}}{\p x_2}\res\vec{e}_{2,k'}\rg)^2\rightharpoonup\Delta\beta_\infty+\lf(\frac{\p\vec{n}_{\vec{\xi}_{\infty}}}{\p x_1}\res\vec{e}_{1,\infty}\rg)^2+\lf(\frac{\p\vec{n}_{\vec{\xi}_{\infty}}}{\p x_2}\res\vec{e}_{2,\infty}\rg)^2\\[5mm]
\ds\frac{\p^2\beta_\infty}{\p x_1\p x_2}=0\quad\quad\mbox{ in }{\mathcal D}'(D^2).
\end{array}
\rg.
\ee
Or in other words, since (\ref{II.11}) holds,
\be
\label{III.37z}
\lf\{
\begin{array}{l}
\ds |\nabla\vec{n}_{\vec{\xi}_{k'}}|^2\ dx_1\,dx_2\rightharpoonup |\nabla\vec{n}_{\vec{\xi}_{\infty}}|^2\ dx_1\,dx_2+\Delta\beta_\infty\ dx_1\,dx_2\\[5mm]
\ds\frac{\p^2\beta_\infty}{\p x_1\p x_2}=0\quad\quad\mbox{ in }{\mathcal D}'(D^2).
\end{array}
\rg.
\ee
The defect measure is then given by the laplacian of an $W^{1,2}(D^2)$ function $\beta_\infty$ whose distributional cross derivative ${\p^2\beta_\infty}/{\p x_1\p x_2}$ is zero.
This implies (\ref{I.22}) and theorem~\ref{th-I.3} is proved.\hfill $\Box$

{99}
\end{document}